\documentclass[onecolumn,final,a4paper]{elsarticle}

\usepackage[ruled,vlined,linesnumbered]{algorithm2e}
\usepackage{algorithmic}
\usepackage{amsfonts}
\usepackage{amsmath}
\usepackage{amsopn}
\usepackage{amssymb}
\usepackage{amsthm}
\usepackage{bbm}
\usepackage{booktabs}
\usepackage{color}
\usepackage[useregional]{datetime2} 
\usepackage{epsfig}
\usepackage{fancyhdr}
\usepackage{float}
\usepackage{graphicx}
\usepackage[colorlinks=true,breaklinks=true,linkcolor=lightblue,citecolor=lightblue]{hyperref}
\usepackage[utf8]{inputenc}
\usepackage{latexsym}
\usepackage{lipsum}
\usepackage{makecell}
\usepackage{mathrsfs}
\usepackage{mathtools}
\usepackage{multirow}
\usepackage{rotating}
\usepackage{showkeys}
\usepackage{scalerel,stackengine}
\stackMath
\newcommand\reallywidehat[1]{%
\savestack{\tmpbox}{\stretchto{%
  \scaleto{%
    \scalerel*[\widthof{\ensuremath{#1}}]{\kern-.6pt\bigwedge\kern-.6pt}%
    {\rule[-\textheight/2]{1ex}{\textheight}}
  }{\textheight}%
}{0.5ex}}%
\stackon[1pt]{#1}{\tmpbox}%
}
\usepackage{stmaryrd}
\usepackage{subfigure}
\usepackage{todonotes}
\usepackage{wrapfig}

\let\oldnl\nl
\numberwithin{equation}{section}
\setuptodonotes{fancyline, color=green!30}

\newtheorem{remark}{Remark}[section]

\def\bB{\mathbf{B}}

\def\bE{\mathbf{E}}

\newcommand\bn{\mathbf{n}}

\def\bU{\mathbf{U}}

\def\bv{\mathbf{v}}
\def\bV{\mathbf{V}}

\def\bx{\mathbf{x}}

\newcommand{\HAT}[1]{\widetilde{}}

\def\J{\mathbf{J}}

\newcommand{\nonl}{\renewcommand{\nl}{\let\nl\oldnl}}%

\newcommand{\vertiii}[1]{{\left\vert\kern-0.25ex\left\vert\kern-0.25ex\left\vert #1
    \right\vert\kern-0.25ex\right\vert\kern-0.25ex\right\vert}}

\graphicspath{{./figs/landau/}{./figs/two_stream/}{./figs/time_reversal/}}

\begin{document}

\begin{frontmatter}

\title{A Tensor-Train Discontinuous Galerkin Method \\ for the Vlasov-Maxwell System}
\author[TDIV]{Rujeko Chinomona\corref{cor1}}
\author[TDIV,IIT]{Dibyendu Adak}
\author[TDIV]{William J. Barham}
\author[TDIV]{Duc P. Truong}
\author[SNL]{Nathan V. Roberts}
\author[TDIV]{Kim \O. Rasmussen}
\author[TDIV,SLIC]{Boian S. Alexandrov}

\address[TDIV]{Theoretical Division,
Los Alamos National Laboratory, Los Alamos, NM 87545, USA}
\address[SNL]{Sandia National Laboratories, Albuquerque, NM, USA}
\address[SLIC]{SLIC.AI, Santa Fe, NM, USA}
\address[IIT]{Indian Institute of Technology-Kharagpur, WB 721302, India}

\cortext[cor1]{Corresponding author. Email: crujeko@lanl.gov}
\begin{abstract}

We present a tensor-train discontinuous Galerkin (TT-DG) formulation for the Vlasov--Maxwell system that combines a modal DG discretization with low-rank tensor representations of the phase-space solution and discrete operators. The formulation exploits the tensor-product structure of the DG discretization to perform quadrature, differentiation, nonlinear upwind flux evaluation, and time integration directly in compressed form.

The method is evaluated on several standard 1D2V Vlasov--Maxwell benchmark problems, including the streaming Weibel instability, weak Landau damping, and two-stream instability problems. Across these problems, the TT formulation reproduces the accuracy and conservation behavior of the underlying full-grid DG discretization while substantially reducing memory usage and runtime. For weakly nonlinear problems, compression ratios exceeding $10^4$ are obtained together with significant speedups relative to the full-grid solver. For the strongly nonlinear two-stream instability problem, the TT formulation remains effective despite reduced compressibility caused by fine-scale phase-space filamentation.

These results demonstrate that tensor-train representations provide an effective approach for reducing the computational cost of deterministic DG-based kinetic plasma simulations while retaining the favorable numerical properties of the underlying discretization.

\end{abstract}

\end{frontmatter}


\section{Introduction}

The Vlasov--Maxwell (VM) system is a fundamental kinetic model for collisionless plasmas
describing the coupled evolution of a particle distribution function in phase-space and
electromagnetic fields. In the VM system, the Vlasov equation governs the evolution of the
particle distribution function under electromagnetic forces, while Maxwell’s equations determine
the electric and magnetic fields generated by the charged particles. These equations arise in applications including magnetic confinement fusion
\cite{freidberg2007plasma}, laser--plasma interaction
\cite{gibbon2005short}, beam physics \cite{reiser2008theory},
and space and astrophysical plasmas \cite{boyd2003physics}.

Unlike fluid or magnetohydrodynamic models \cite{hazeltine2004framework,nicholson1983plasma}, the VM system resolves velocity-space effects
such as Landau damping, particle trapping, filamentation, and anisotropy. However, this kinetic
description comes with substantial computational cost. Deterministic numerical solution of the VM
system requires discretizing a phase-space of dimension $d=d_x+d_v$. In the present work, we
focus on the 1D2V setting, where $d_x=1$ and $d_v=2$, giving a three-dimensional phase-space.
Even in this reduced setting, storage and computational cost scale as $\mathcal{O}(N^3)$ on a
uniform grid with $N$ points per dimension, and the problem rapidly becomes more expensive
under refinement. Extending to higher-dimensional VM systems, where $d$ reaches $4$, $5$, or
$6$, makes full-grid simulation prohibitively expensive in both memory and runtime. Developing
numerical methods that retain the accuracy and conservation properties of the VM system while
reducing this computational cost remains a major challenge in computational plasma physics.

We consider the VM system for a single species of nonrelativistic electrons, 
while the ions are treated as a fixed, spatially uniform background. 
Using a common nondimensionalization \cite{li2019solving}, 
with time scaled by the inverse plasma frequency, velocity by the speed of
light, and length by the electron skin depth, the dimensionless VM equations are
given by
\begin{align}
\partial_t f + v \cdot \nabla_x f + (E + v \times B)\cdot \nabla_v f &= 0,
\label{eq:vlasov} \\
\partial_t E &= \nabla_x \times B - J,
\label{eq:maxwellE} \\
\partial_t B &= -\nabla_x \times E,
\label{eq:maxwellB} \\
\nabla_x \cdot E &= \rho - \rho_i,
\qquad
\nabla_x \cdot B = 0.
\label{eq:constraints}
\end{align}
The density and current density are defined by
\begin{align}
\rho(x,t) = \int_{\Omega_v} f(x,v,t)\,dv, \qquad
J(x,t) = \int_{\Omega_v} f(x,v,t)\,v\,dv ,
\end{align}
while $\rho_i$ is the ion density. Here, $f(x,v,t)$ denotes the probability distribution function for finding an electron at position
$x$ with velocity $v$ at time $t$. The electric and magnetic fields are $E$ and $B$ respectively. The spatial domain is denoted by $\Omega_x$, while $\Omega_v$
denotes the velocity domain. The VM system conserves several important physical quantities,
including the total particle number
\[
\int_{\Omega_x}\int_{\Omega_v} f\,dv\,dx,
\]
and the total energy
\[
\mathcal{E}
=
\frac12
\int_{\Omega_x}\int_{\Omega_v}
f|v|^2\,dv\,dx
+
\frac12
\int_{\Omega_x}
\left(|E|^2+|B|^2\right)\,dx,
\]
which consists of the kinetic and electromagnetic energy contributions.

\begin{remark}
The skin-depth/light-speed normalization is most natural for electromagnetic examples, where
transverse Maxwell dynamics are part of the leading-order model and nonrelativistic distributions
satisfy $|v|\ll 1$. Electrostatic tests such as weak Landau damping are usually interpreted instead
in Debye-length and thermal-speed units. In that scaling the small parameter
$\varepsilon=v_T/c$ appears in the full VM system; 
in particular, Amp\`ere's law carries the factor $\varepsilon^{-2}$ on the Maxwell curl term \cite{crouseilles2016asymptotic}. For data lying on the electrostatic invariant
submanifold---that is, data with longitudinal electric field, zero magnetic field, and no transverse
current under the evolution---the transverse Maxwell dynamics are not excited at the continuum
level. This is why full VM codes can run electrostatic benchmarks without resolving the transverse
light-wave scale, provided the discretization preserves this invariance. In
Section~\ref{sec:sec_num_exp}, we interpret the electromagnetic tests using the
skin-depth/light-speed normalization and the electrostatic weak Landau damping test using the
Debye-length/thermal-speed normalization.
\end{remark}

A wide variety of numerical methods have been developed for the Vlasov-Maxwell system,
including particle-in-cell methods~\cite{kraus2017gempic,xiao2018structure},
semi-Lagrangian methods~\cite{sonnendrucker1999semi,Bostan2009convergence},
spectral methods~\cite{holloway1996spectral,delzanno2015multidimensional},
finite difference methods~\cite{sircombe2009valis, shiroto2019quadratic}, phase-space finite volume methods~\cite{filbet2001conservative},
and discontinuous Galerkin (DG) methods~\cite{cheng2014discontinuous,cheng2014energy}.
Among deterministic approaches, DG methods are particularly attractive because they naturally
support high-order accuracy, local conservation, geometric flexibility, and structure-preserving
formulations. Critically for the present work, modal DG discretizations on tensor-product meshes
produce solution coefficients that organize naturally into multidimensional arrays, providing a
direct foundation for low-rank tensor compression.

Tensor-network (TN) methods are structured low-rank representations designed to mitigate the curse of
dimensionality by factorizing high-dimensional tensors into networks of lower-order tensors.
Among these, the tensor-train (TT) format introduced by Oseledets and
Tyrtyshnikov~\cite{oseledets2011tt,oseledets2010ttcross} is particularly well suited for
tensor-product DG discretizations. In the TT representation, a $d$-dimensional tensor is
decomposed into a sequence of coupled three-dimensional cores whose interconnection ranks $r$
determine the approximation accuracy. When these ranks remain moderate, storage cost reduces
from $\mathcal{O}(N^d)$ to $\mathcal{O}(dNr^2)$, reducing the exponential scaling of full-grid
representations to linear scaling in dimension.

Low-rank approaches for Vlasov-type equations have received growing attention,
including TT semi-Lagrangian
methods~\cite{kormann2015semilagrangian}, adaptive tensor
decompositions~\cite{EHRLACHER2017285}, dynamical low-rank
projector-splitting methods~\cite{einkemmer2018lowrank}, low-rank tensor
formulations with hierarchical Tucker extensions for high-dimensional
problems~\cite{guo2022lowrank,guo2024conservative,guo2024lomacdg}, and
adaptive-rank semi-Lagrangian
methods~\cite{zheng_semi-lagrangian_2025}. 
Related TN approaches have also been developed for other high-dimensional
transport equations, including TT/QTT formulations of the Boltzmann neutron transport
equation in which both the angular flux and discrete transport operators are represented in TT or QTT format~\cite{truong2024NTE,ortega2026nuclear}.
More recently, tensor-network
formulations have begun to emerge for the Vlasov--Maxwell
system~\cite{ye2024quantized,TegbyOskar2023AHTS}. However, a TT formulation
combined with a DG discretization for the full Vlasov--Maxwell
system, which requires handling variable-coefficient electromagnetic
transport terms and nonlinear flux evaluations within the TT
framework, has not previously been developed.

To the best of our knowledge, this work presents the first tensor-train discontinuous Galerkin
(TT-DG) formulation for the Vlasov--Maxwell system. We reformulate the modal DG discretization
to expose the underlying tensor-product structure of the phase-space solution and discrete operators,
develop TT realizations of the DG volume and surface operators including the variable-coefficient
Lorentz force terms and nonlinear upwind fluxes, and investigate the influence of tensor layout
and phase-space ordering on compression efficiency and runtime. Through benchmark studies of
streaming Weibel instability, weak Landau damping, and two-stream instability, we demonstrate
that the proposed formulation achieves compression ratios up to $\mathcal{O}(10^4)$ and speedups
up to $\mathcal{O}(10^2)$ relative to the full-grid DG solver while preserving the accuracy and
conservation properties of the underlying discretization.

The remainder of this paper is organized as follows. In Section~\ref{sec:sec2}, we present the DG
discretization of the Vlasov--Maxwell system and introduce the tensor-product modal formulation
used throughout the paper. Section~\ref{sec:sec_tt} develops the tensor-train realization of the DG
discretization, including TT representations of the solution coefficients, operator evaluations,
flux computations, and time integration. Section~\ref{sec:sec_num_exp} presents numerical experiments for several
benchmark problems, including streaming Weibel instability, weak Landau damping, and
two-stream instability, together with studies of TT compression, tensor ordering, runtime, and
conservation behavior. Finally, Section~\ref{sec:sec_conclusion} summarizes the main conclusions and discusses future
research directions.
\section{Discontinuous Galerkin discretization and tensor-product structure} \label{sec:sec2}
We present a modal discontinuous Galerkin (DG) discretization of the Vlasov--Maxwell system, coupled with an explicit time integrator following \cite{cheng2014discontinuous,cheng2014energy}. After introducing the notation and discrete spaces, we derive the semi-discrete DG formulation of \eqref{eq:vlasov}--\eqref{eq:maxwellB}, and subsequently construct the fully discrete scheme.
The divergence constraints \eqref{eq:constraints} are not enforced explicitly.

\subsection{Notation}\label{subsec:sec2_notation}

Let $d_x$ and $d_v$ denote the number of spatial and velocity dimensions, respectively. We first define the computational mesh. The spatial and velocity domains are decomposed coordinate-wise as
\[
\Omega_x = \Omega_{x_1} \times \cdots \times \Omega_{x_{d_x}},
\qquad
\Omega_v = \Omega_{v_1} \times \cdots \times \Omega_{v_{d_v}}.
\]
For each $\ell=1,\dots,d_x$, let
\[
\mathcal{T}_h^{x_\ell} = \{ K_{x_\ell,i_\ell} \}_{i_\ell=1}^{N_{x_\ell}}
\]
be a partition of $\Omega_{x_\ell}$, and for each $m=1,\dots,d_v$, let
\[
\mathcal{T}_h^{v_m} = \{ K_{v_m,j_m} \}_{j_m=1}^{N_{v_m}}
\]
be a partition of $\Omega_{v_m}$. Here $h$ denotes the maximum interval length of all one-dimensional spatial and velocity elements in the mesh. The multi-dimensional elements are given by
\[
K_x = \prod_{\ell=1}^{d_x} K_{x_\ell,i_\ell},
\qquad
K_v = \prod_{m=1}^{d_v} K_{v_m,j_m},
\]
with multi-indices $\mathbf{i}=(i_1,\dots,i_{d_x})$ and $\mathbf{j}=(j_1,\dots,j_{d_v})$. The phase-space elements are
\[
K = K_x \times K_v,
\qquad
\mathcal{T}_h = \{ K : \mathbf{i}, \mathbf{j} \}.
\]

For completeness, we briefly describe the edge and trace notation used in the DG formulation. Edges arise as interfaces between neighboring phase-space elements in each coordinate direction, induced by the one-dimensional meshes $\mathcal{T}_h^{x_\ell}$ and $\mathcal{T}_h^{v_m}$. Let $e = K^+ \cap K^-$ denote an interior edge shared by two elements $K^\pm$. We denote the traces of a scalar function $g$ on $e$ by $g^\pm = g|_{K^\pm}$ and the corresponding outward unit normals by $\bn^\pm$.
These traces are used to define averages and jumps. For a scalar function $g$ and a vector-valued function $\bU$, we set
\[
\{g\} = \tfrac{1}{2}(g^+ + g^-),
\qquad
[g] = g^+ \bn^+ + g^- \bn^-,
\]
\[
\{\bU\} = \tfrac{1}{2}(\bU^+ + \bU^-),
\qquad
[\bU] = \bU^+ \cdot \bn^+ + \bU^- \cdot \bn^-.
\]

We next introduce the discrete spaces. Let $p \ge 0$ be a fixed polynomial degree, and for any interval $I \subset \mathbb{R}$, let $P^p(I)$ denote the space of polynomials of degree at most $p$ on $I$. For a spatial element $K_x$, define the tensor-product space
\[
Q^p(K_x) = \bigotimes_{\ell=1}^{d_x} P^p(K_{x_\ell,i_\ell}),
\]
and for a phase-space element $K$,
\[
Q^p(K)
=
\bigotimes_{\ell=1}^{d_x} P^p(K_{x_\ell,i_\ell})
\otimes
\bigotimes_{m=1}^{d_v} P^p(K_{v_m,j_m}).
\]
Here $Q^p(K)$ denotes the tensor-product polynomial space of degree at most $p$ in each coordinate direction.

The DG phase-space is
\[
V_h^p = \{ g \in L^2(\Omega) : g|_K \in Q^p(K), \ \forall K \in \mathcal{T}_h \},
\]
and the vector-valued physical-space is
\[
X_h^p = \{ \bU \in [L^2(\Omega_x)]^{d_x} : \bU|_{K_x} \in [Q^p(K_x)]^{d_x}, \ \forall K_x \in \mathcal{T}_h^x \}.
\]

The local spaces are represented using a tensor-product modal basis constructed from mapped Legendre polynomials. Let $\{P_i(\xi)\}_{i=0}^p$ denote the Legendre polynomials on $[-1,1]$, and let
\[
\xi_{x_\ell}(x_\ell) = \frac{2(x_\ell - x_{\ell,L})}{\Delta x_\ell} - 1
\]
be the affine mapping from $K_{x_\ell,i_\ell}$ to the reference interval. We define
\[
\phi^{x_\ell}_{i}(x_\ell) = P_i(\xi_{x_\ell}(x_\ell)), 
\qquad
\phi^{v_m}_{j}(v_m) = P_j(\xi_{v_m}(v_m)).
\]
With this choice, $Q^p(K)$ admits the representation
\[
Q^p(K)
=
\operatorname{span}
\left\{
\prod_{\ell=1}^{d_x} \phi^{x_\ell}_{i_\ell}(x_\ell)
\prod_{m=1}^{d_v} \phi^{v_m}_{j_m}(v_m)
\right\}.
\]
The multidimensional basis is formed as a tensor product of one-dimensional basis functions, the DG expansion coefficients are naturally organized as a multidimensional tensor. This tensor-product structure is particularly well suited for tensor-train compression and tensor-product operator evaluation. In addition, the $L^2$-orthogonality of the Legendre basis yields diagonal element mass matrices and simplifies the resulting operator structure.

\subsection{Semi-discrete DG formulation} \label{subsec:sec2_semi_discrete}

We derive the semi-discrete DG formulation of the Vlasov--Maxwell system. We seek DG approximations
\[
f_h(t) \in V_h^p, 
\qquad 
\bE_h(t), \bB_h(t) \in X_h^p,
\]
for all $t>0$. The tensor-product modal basis introduced in Section 2.1 induces a multidimensional coefficient structure for the DG solution, which later forms the basis for the tensor-train representation. On each phase-space element $K \in \mathcal{T}_h$, the numerical distribution function satisfies
\[
f_h|_K(\bx,\bv,t) \in Q^p(K),
\]
and admits the modal expansion
\begin{equation}\label{eq:modal_expansion}
f_h|_K(\bx,\bv,t)
=
\sum_{\mathbf{i},\mathbf{j}}
\reallywidehat{f}^{K}_{\mathbf{i},\mathbf{j}}(t)
\prod_{\ell=1}^{d_x}
\phi^{x_\ell}_{i_\ell}(x_\ell)
\prod_{m=1}^{d_v}
\phi^{v_m}_{j_m}(v_m),
\end{equation}
where $\reallywidehat{f}^{K}_{\mathbf{i},\mathbf{j}}(t)$ are the time-dependent modal DG coefficients associated with the tensor-product basis functions on the element $K$, $\mathbf{i}=(i_1,\dots,i_{d_x})$, $\mathbf{j}=(j_1,\dots,j_{d_v})$, and $0\le i_\ell,j_m\le p$.
The electromagnetic fields $\bE_h$ and $\bB_h$ are vector-valued DG functions defined only on the physical-space mesh, reflecting their dependence on the spatial coordinates $\bx$ alone. 
 
The semi-discrete formulation is obtained by multiplying each equation of the Vlasov–Maxwell system by appropriate test functions $\psi \in V_h^p$, integrating over elements, and applying integration by parts to expose interface flux contributions. Because the DG approximation is generally discontinuous across element interfaces, we introduce single-valued numerical fluxes on interfaces, denoted throughout by an over-hat notation. On spatial and velocity faces, we write the corresponding outward normals as $\bn_x$ and $\bn_v$. For the Vlasov equation, this yields, for all $\psi \in V_h^p$,
\begin{align}
\int_K \partial_t f_h \, \psi \, d\bx d\bv
&- \int_K f_h \, \bv \cdot \nabla_x \psi \, d\bx d\bv
- \int_K f_h \, (\bE_h + \bv \times \bB_h) \cdot \nabla_{\bv} \psi \, d\bx d\bv \\
&+ \int_{\partial K_x \times K_v} \reallywidehat{f_h \, \bv} \cdot \bn_x \, \psi \, ds_x d\bv
+ \int_{K_x \times \partial K_v} \reallywidehat{f_h (\bE_h + \bv \times \bB_h)} \cdot \bn_v \, \psi \, d\bx ds_v
= 0. \nonumber
\end{align}
The first two volume integrals correspond to spatial and velocity-space transport, while the surface terms introduce numerical fluxes that weakly couple neighboring elements across interfaces.

Applying the same procedure to Maxwell's equations, we obtain, for all $\bU,\bV \in X_h^p$,
\begin{align}
\int_{K_x} \partial_t \bB_h \cdot \bU \, d\bx
&= - \int_{K_x} \bE_h \cdot (\nabla_x \times \bU) \, d\bx
+ \int_{\partial K_x} \reallywidehat{\bn \times \bE_h} \cdot \bU \, ds, \\
\int_{K_x} \partial_t \bE_h \cdot \bV \, d\bx
&= \int_{K_x} \bB_h \cdot (\nabla_x \times \bV) \, d\bx
- \int_{K_x} \J_h \cdot \bV \, d\bx \nonumber \\
&\quad - \int_{\partial K_x} \reallywidehat{\bn \times \bB_h} \cdot \bV \, ds,
\end{align}
where
\[
\J_h(\bx,t) = \int_{\Omega_v} f_h(\bx,\bv,t)\,\bv\, d\bv.
\]

The numerical fluxes are constructed from the jumps and averages introduced above. For the Vlasov equation, we employ upwind numerical fluxes, which introduce dissipation that improves stability and suppresses nonphysical oscillations associated with phase-space filamentation in kinetic simulations \cite{cheng2014energy}. The resulting spatial and velocity fluxes are
\begin{equation}
\reallywidehat{f_h \bv}
=
\{f_h \bv\}
+
\frac{|\bv\cdot\bn_x|}{2}[f_h],
\qquad
\reallywidehat{f_h (\bE_h+\bv\times\bB_h)}
=
\{f_h(\bE_h+\bv\times\bB_h)\}
+
\frac{|(\bE_h+\bv\times\bB_h)\cdot\bn_v|}{2}[f_h].
\end{equation}
For Maxwell's equations, we employ the alternating numerical fluxes
\begin{equation}
\reallywidehat{\bn \times \bE_h} = \bn \times \bE_h^{+}, 
\qquad
\reallywidehat{\bn \times \bB_h} = \bn \times \bB_h^{-},
\end{equation}
or equivalently
\begin{equation}
\reallywidehat{\bn \times \bE_h} = \bn \times \bE_h^{-}, 
\qquad
\reallywidehat{\bn \times \bB_h} = \bn \times \bB_h^{+}.
\end{equation}
An alternative choice is the central flux. However, as shown in \cite{cheng2014energy}, the alternating flux provides improved numerical accuracy for the VM system while retaining the desired conservation properties, and is therefore adopted in this work.

For later use in the TT formulation, we introduce compact operator notation for the semi-discrete DG system. For $f_h,\psi \in V_h^p$, define
\[
\mathcal{M}_f(\partial_t f_h,\psi)
=
\sum_{K\in\mathcal{T}_h}
\int_K \partial_t f_h \, \psi \, d\bx d\bv ,
\]
and let $\mathcal{A}_h(f_h;\bE_h,\bB_h,\psi)$ denote the sum of the corresponding volume and surface contributions.

For $\bE_h,\bB_h \in X_h^p$, we define the mass forms $\mathcal{M}_E$, $\mathcal{M}_B$, and the discrete curl operator
\[
\mathcal{C}_h(\bE_h,\bU)
=
\sum_{K_x\in\mathcal{T}_h^x}
\int_{K_x} \bE_h \cdot (\nabla_x \times \bU)\, d\bx
-
\sum_{K_x\in\mathcal{T}_h^x}
\int_{\partial K_x}
\reallywidehat{\bn \times \bE_h} \cdot \bU \, ds,
\]
and analogously for $\mathcal{C}_h(\bB_h,\bV)$.

We also define the $L^2$ inner product on $\Omega_x$ by
\[
(\J_h,\bV)_{\Omega_x}
=
\int_{\Omega_x} \J_h(\bx,t)\cdot \bV(\bx)\, d\bx.
\]

With this notation, the global semi-discrete DG formulation reads: find
\[
f_h(t)\in V_h^p,
\qquad
\bE_h(t),\bB_h(t)\in X_h^p,
\]
such that for all $\psi\in V_h^p$ and $\bU,\bV\in X_h^p$,
\begin{equation}\label{eq:semi_discrete}
\begin{aligned}
&\mathcal{M}_f(\partial_t f_h,\psi)
+
\mathcal{A}_h(f_h;\bE_h,\bB_h,\psi)
=0, \\
&\mathcal{M}_B(\partial_t\bB_h,\bU)
+
\mathcal{C}_h(\bE_h,\bU)
=0, \\
&\mathcal{M}_E(\partial_t\bE_h,\bV)
-
\mathcal{C}_h(\bB_h,\bV)
+
(\J_h,\bV)_{\Omega_x}
=0.
\end{aligned}
\end{equation}

\subsection{Fully discrete scheme} \label{subsec:sec2_fully_discrete}

For our temporal discretization, we employ the second-order explicit Scheme-1 from Cheng et al. \cite{cheng2014energy}, which combines a leapfrog discretization for Maxwell’s equations with a second-order explicit Runge--Kutta discretization for the Vlasov equation. The method is designed to respect the Hamiltonian structure of the Vlasov--Maxwell system and exhibits favorable long-time energy behavior. While the scheme does not exactly conserve the total energy, it exactly preserves a modified discrete energy that remains close to the physical total energy \cite{cheng2014energy}.

The staggered structure of the scheme separates the Vlasov and Maxwell updates while preserving the explicit coupling through the midpoint current density.
Let $\Delta t>0$ be the time step and define the discrete times $t^n = n\Delta t$ for integers $n \ge 0$. Given $f_h^n$, $\bE_h^n$, and $\bB_h^n$, we compute the intermediate states $f_h^{n+\frac12}$ and $\bB_h^{n+\frac12}$, followed by the updated solutions $\bE_h^{n+1}$, $\bB_h^{n+1}$, and $f_h^{n+1}$ according to
\begin{subequations}\label{eq:scheme1}
\begin{align}
&\mathcal{M}_f
\left(
\frac{f_h^{n+\frac12}-f_h^n}{\Delta t/2},\psi
\right)
+
\mathcal{A}_h(f_h^n;\bE_h^n,\bB_h^n,\psi)
=0,
\label{eq:scheme1a}
\\
&\mathcal{M}_B
\left(
\frac{\bB_h^{n+\frac12}-\bB_h^n}{\Delta t/2},\bU
\right)
+
\mathcal{C}_h(\bE_h^n,\bU)
=0,
\label{eq:scheme1b}
\\
&\mathcal{M}_E
\left(
\frac{\bE_h^{n+1}-\bE_h^n}{\Delta t},\bV
\right)
-
\mathcal{C}_h(\bB_h^{n+\frac12},\bV)
+
(\J_h^{n+\frac12},\bV)_{\Omega_x}
=0,
\label{eq:scheme1c}
\\
&\mathcal{M}_B
\left(
\frac{\bB_h^{n+1}-\bB_h^{n+\frac12}}{\Delta t/2},\bU
\right)
+
\mathcal{C}_h(\bE_h^{n+1},\bU)
=0,
\label{eq:scheme1d}
\\
&\mathcal{M}_f
\left(
\frac{f_h^{n+1}-f_h^n}{\Delta t},\psi
\right)
+
\mathcal{A}_h
\left(
f_h^{n+\frac12};
\frac{\bE_h^n+\bE_h^{n+1}}{2},
\bB_h^{n+\frac12},
\psi
\right)
=0,
\label{eq:scheme1e}
\end{align}
\end{subequations}
for all $\psi \in V_h^p$ and $\bU,\bV \in X_h^p$, where
\[
\J_h^{n+\frac12}(\bx)
=
\int_{\Omega_v} f_h^{n+\frac12}(\bx,\bv)\,\bv\,d\bv.
\]

This formulation is particularly convenient for the TT implementation developed later, since the Vlasov updates can be performed directly on the compressed coefficient tensors while the Maxwell fields remain in their standard DG representation.
\subsection{Tensor-product representation of DG operators}
\label{subsec:sec2_tensor_product}

To prepare for the tensor-train realization developed in
Section~\ref{sec:sec_tt}, we briefly describe the tensor-product
operator structure underlying the full-grid DG implementation.
These implementation-oriented details do not modify the DG formulation,
but expose the separable algebraic structure used in the TT
representation.

The tensor-product modal basis introduced in
Section~\ref{subsec:sec2_semi_discrete} induces a corresponding
tensor-product structure in the discrete DG operators. Let
\[
\mathcal Z
=
\{x_1,\ldots,x_{d_x},v_1,\ldots,v_{d_v}\}
\]
denote the set of phase-space coordinate directions. For each
$z\in\mathcal Z$, let
$\{z_q\}_{q=1}^{n_q}$ and
$\{\omega_q^z\}_{q=1}^{n_q}$
denote the quadrature points and weights on the one-dimensional cell
$K_z$, and let $n_b=p+1$ denote the number of one-dimensional modal
basis functions.
We define the basis-evaluation, derivative-evaluation, quadrature-weight,
and mass matrices by
\[
\Phi_z \in \mathbb R^{n_q\times n_b},
\qquad
D_z \in \mathbb R^{n_q\times n_b},
\qquad
W_z \in \mathbb R^{n_q\times n_q},
\qquad
M_z \in \mathbb R^{n_b\times n_b},
\]
with entries
\[
(\Phi_z)_{qr} = \phi_r^z(z_q),
\qquad
(D_z)_{qr} = \partial_z \phi_r^z(z_q),
\qquad
W_z = \operatorname{diag}(\omega_1^z,\ldots,\omega_{n_q}^z),
\]
and
\[
(M_z)_{rs}
=
\int_{K_z}
\phi_r^z(z)\phi_s^z(z)\,dz
\approx
(\Phi_z^\top W_z \Phi_z)_{rs}
\]
respectively.
The tensor-product basis implies that multidimensional DG operators
factor into tensor (Kronecker) products of one-dimensional operators.
In particular, the element mass matrix associated with the Vlasov
equation factors as
\[
M_K^f
=
\bigotimes_{z\in\mathcal Z} M_z,
\]
where $\otimes$ denotes the Kronecker product over the phase-space
directions in $\mathcal Z$. Consequently, the mass form
$\mathcal M_f(\partial_t f_h,\psi)$ introduced in
Section~\ref{subsec:sec2_semi_discrete} admits a separable tensor-product
representation in terms of one-dimensional operators.
The same tensor-product structure applies to the DG volume and surface
operators appearing in
$\mathcal A_h(f_h;\bE_h,\bB_h,\psi)$.
In particular, multidimensional quadrature, differentiation, and face
evaluations can be represented through tensor-product contractions
involving the one-dimensional matrices
$\Phi_z$, $D_z$, $W_z$, together with the corresponding
face-evaluation operators.
For example, if
$\reallywidehat{\mathbf f}^{\,K}$ denotes the vector of DG modal coefficients
associated with the expansion \eqref{eq:modal_expansion} on an element
$K$, then the corresponding values of $f_h$ at the tensor-product
quadrature points are given by
\[
\mathbf f_q^K
=
\left(
\bigotimes_{z\in\mathcal Z}\Phi_z
\right)
\reallywidehat{\mathbf f}^{\,K}.
\]
Similarly, directional differentiation in a coordinate direction $z$
corresponds to replacing the factor $\Phi_z$ by $D_z$ in the associated
tensor-product operator. 
This separable operator structure is the key property that enables efficient low-rank tensor representations of the DG discretization. Exposing tensor-product structure at the discretization level allows quadrature, differentiation, and inverse mass operations to be applied one coordinate direction at a time, forming the basis for the TT implementation developed in Section~\ref{sec:sec_tt}. For Maxwell's equations, the fields depend only on the physical-space variables and therefore admit an analogous tensor-product structure over the spatial directions.

Exploiting such tensor-product operators is a common theme in tensor-network PDE solvers. For example, TT formulations for three-dimensional Maxwell wave propagation leverage separable spatial operators~\cite{manzini2023Maxwell}, while TT isogeometric analysis employs tensor-product spline spaces together with compressed geometry-dependent coefficients to assemble PDE operators on nontrivial geometries~\cite{tran2026TTIGA}.

\section{Tensor-train formulation}\label{sec:sec_tt}
The tensor-product DG discretization developed in Section \ref{sec:sec2} induces a separable multidimensional coefficient structure for both the phase-space solution and the associated discrete operators. In this section, we exploit this structure by representing the DG coefficient tensor in tensor-train (TT) format and performing the DG operator evaluations directly in compressed form. This enables low-rank realization of the Vlasov equation without explicitly forming full multidimensional arrays, substantially reducing the storage and computational cost of the kinetic solve.

We focus on the one-dimensional in space, two-dimensional in velocity (1D2V) Vlasov--Maxwell system with phase-space variables $(x_2,v_{x1},v_{x2})$. Since the electromagnetic fields depend only on the spatial  coordinate, the TT representation is applied only to the phase-space discretization of the Vlasov equation, while Maxwell's equations retain the standard full-grid DG representation introduced in Section \ref{sec:sec2}.

The DG formulation requires evaluating transport operators and nonlinear upwind fluxes directly within the TT framework while maintaining a compressed representation throughout the computation. This introduces several additional considerations, including tensor ordering, quadrature evaluation in TT form, nonlinear flux construction, TT-cross approximation, and low-rank time integration. These components are developed in the following subsections.

\subsection{Tensor-train representation and ordering}
\label{subsec:tt-representation-ordering}
The tensor-product DG discretization organizes the modal
coefficients of $f_h$ into a multidimensional array with one element
index and one local basis index per phase-space direction. In the
1D2V setting with phase-space variables
$(x_2,v_{x_1},v_{x_2})$, the modal expansion
\eqref{eq:modal_expansion} therefore produces a six-dimensional
coefficient tensor. Each phase-space direction contributes both an element index and a local modal-basis index, producing two tensor modes per coordinate direction. Since each one-dimensional polynomial space has
$n_b=p+1$ modal basis functions, the tensor dimensions consist of the
element counts $(N_{x_2},N_{v_{x_1}},N_{v_{x_2}})$ together with the
corresponding modal dimensions $(n_b,n_b,n_b)$.  
We denote the resulting full-grid coefficient tensor by $\mathcal{F}$ and its compressed TT-format approximation $\mathcal{F}^{\mathrm{TT}}$. The TT approximation is constructed to satisfy a prescribed truncation tolerance $\varepsilon_{\mathrm{TT}}$, which controls the accuracy of low-rank compression and recompression operations throughout the computation. The tensor-train representation~\cite{oseledets2011tt} is
\[
\mathcal{F}^{\mathrm{TT}}(i_1,\dots,i_6)
=
\sum_{\alpha_1,\dots,\alpha_5}
G_1(1,i_1,\alpha_1)
G_2(\alpha_1,i_2,\alpha_2)\cdots
G_6(\alpha_5,i_6,1),
\]
where \(G_k\) are the TT cores and \(\alpha_k\) index the intermediate TT ranks.
Different orderings expose different correlation structures to the TT compression algorithm and can therefore substantially affect both achievable ranks and computational cost.

We consider two natural layouts. Let $\pi=(\pi_1,\pi_2,\pi_3)$ be a permutation of $(x_2,v_{x_1},v_{x_2})$. The interleaved ordering is
\begin{equation}
[\;N_{\pi_1},\, n_b,\, N_{\pi_2},\, n_b,\, N_{\pi_3},\, n_b\;],
\label{eq:tt-layout-interleaved}
\end{equation}
while the grouped ordering is
\begin{equation}
[\;N_{\pi_1},\, N_{\pi_2},\, N_{\pi_3},\, n_b,\, n_b,\, n_b\;].
\label{eq:tt-layout-grouped}
\end{equation}
For example, the interleaved layout with
\[
\pi=(v_{x_1},x_2,v_{x_2})
\]
corresponds to the ordering
\[
[N_{v_{x_1}},\, n_b,\,
  N_{x_2},\, n_b,\,
  N_{v_{x_2}},\, n_b].
\]
In the interleaved layout, element and modal indices associated with the same coordinate direction are placed adjacently, whereas the grouped layout separates all element indices from all modal indices.
The impact of layout and permutation choice on compression and computational cost is examined in
Section~\ref{subsec:tt_layout_study}.

\subsection{TT-cross interpolation for nonlinear tensor functions}
\label{subsec:tt-cross}

In addition to standard TT arithmetic, the DG discretization requires entrywise nonlinear operations on tensors represented in TT format. A typical example arises in the upwind numerical flux, where the split form
\[
T^\pm = \frac12(T \pm |T|)
\]
requires the entrywise absolute value of the normal transport speed \(T\). 
Since the nonlinear map \(T \mapsto |T|\) is not a standard TT operation, explicitly constructing \(|T|\) would generally require forming the full tensor, which is infeasible in high dimensions. We therefore approximate |T| directly in TT format using TT-cross interpolation, which adaptively samples selected tensor entries without assembling the full array.

Here, TT-cross refers to the DMRG-cross/greedy-cross interpolation strategy used in Savostyanov et al. \cite{Savostyanov:2014} and Dolgov et al. \cite{dolgov2019parallel}, rather than the original one-site TT-cross algorithm of Oseledets and Tyrtyshnikov~\cite{oseledets2010ttcross}. The original TT-cross method uses sequential maxvol-type index selection, typically with prescribed TT ranks. By contrast, the greedy DMRG-cross variant selects new pivots directly from
local two-core superblocks, without performing a maxvol search. These pivots are appended to the interpolation sets, adaptively increasing the TT ranks as needed. This avoids the local cost of maxvol and is useful for nonlinear tensor functions whose effective ranks are not known in advance.

Cross interpolation has become a standard tool for evaluating nonlinear or otherwise nonseparable tensor functions without assembling the full tensor. Beyond the present application, TT-cross has been used for high-dimensional integration~\cite{dolgov2019parallel,alexandrov2023challenging}, nonlinear flux and WENO operator evaluation in tensor-train compressible-flow solvers~\cite{danis2025WENO}, and adaptive evaluation of configurational integrals in statistical mechanics~\cite{truong2025breaking}. These applications demonstrate that cross interpolation provides an efficient and practical mechanism for constructing TT representations of nonlinear tensor functions from adaptively sampled entries.

\subsection{TT evaluation of DG operators}
\label{subsec:tt-dg-operator}
This subsection describes the TT realization of the Vlasov DG right-hand side
associated with the semi-discrete formulation \eqref{eq:semi_discrete}. In the
full-grid formulation, the Vlasov equation has the form
\[
M_f(\partial_t f_h,\psi)
+
\mathcal{A}_h(f_h;\bE_h,\bB_h,\psi)=0,
\]
so the algebraic right-hand side is obtained by assembling the volume and
surface contributions in $\mathcal{A}_h$, applying the appropriate sign, and then applying the inverse mass operator.

The TT formulation replaces the full-grid tensor operations appearing in
$\mathcal{A}_h$ with their TT counterparts. Using the coefficient tensor
ordering introduced in Section~\ref{subsec:tt-representation-ordering}, we
describe the TT realization of the volume operators, surface flux operators,
and inverse mass application that together produce the TT representation of
$\partial_t f_h$.

\paragraph{Volume terms}
Volume contributions are computed by applying the one-dimensional DG operators
$\Phi_z$, $D_z$, and $W_z$ introduced in Section~\ref{subsec:sec2_tensor_product} directly to the TT cores. 
Let $G_z$ denote the modal-basis core associated with direction $z\in\{x_2,v_{x_1},v_{x_2}\}$. 
Evaluating the DG solution at quadrature points amounts to replacing each
modal-basis core $G_z$ by a transformed core $\widetilde G_z$. Using
MATLAB-style indexing notation, the transformation is
\[
\widetilde{G}_z(:,q_z,:)
=
\sum_{\ell_z=1}^{n_b}
G_z(:,\ell_z,:)\Phi_z(q_z,\ell_z),
\]
where $G_z(:,\ell_z,:)$ denotes the matrix obtained by fixing the modal index
to $\ell_z$.
After this update is carried out in all three directions, the resulting TT tensor denoted
$\mathcal{F}^{\mathrm{TT,quad}}$ represents the DG solution sampled on the tensor-product quadrature grid within each phase-space element, i.e. $\mathcal{F}^{\mathrm{TT,quad}}\approx f_h(\bx_q,\bv_q)$.

The directional transport speeds are assembled from one-dimensional
factors on the quadrature grid
\[
T_{x_2}=v_{x_2},\qquad
T_{v_{x_1}}=E_{x_1}+v_{x_2}B_{x_3},\qquad
T_{v_{x_2}}=E_{x_2}-v_{x_1}B_{x_3}.
\]
In the TT representation, these lower-dimensional factors are expanded to the same
tensor-product quadrature layout as $\mathcal{F}^{\mathrm{TT,quad}}$ by
replicating them trivially in the coordinates on which they do not depend.
This yields TT tensors $T_z^{\mathrm{TT}}$ that correspond with
$\mathcal{F}^{\mathrm{TT,quad}}$.
The directional quadrature integrands are then formed by the elementwise (Hadamard) products
\[
\mathcal{G}^{\mathrm{TT}}_{z}
=
T_{z}^{\mathrm{TT}}
\odot
\mathcal{F}^{\mathrm{TT,quad}},
\qquad
z\in\{x_2,v_{x_1},v_{x_2}\}.
\]
Each directional volume contribution is obtained by projecting the corresponding quadrature integrand back to modal form, using the derivative test matrix in the differentiated coordinate and value test matrices in the other two coordinates. For example, the $v_{x_1}$-direction volume term is
\begin{equation}
\mathcal{R}_{\mathrm{vol},v_{x_1}}^{\mathrm{TT}}
=
\left(\Phi_{x_2}^{\top}W_{x_2}\otimes
D_{v_{x_1}}^{\top}W_{v_{x_1}}\otimes
\Phi_{v_{x_2}}^{\top}W_{v_{x_2}}\right)\mathcal{G}^{\mathrm{TT}}_{v_{x_1}},
\end{equation}
with analogous expressions for the $x_2$- and $v_{x_2}$-direction terms. Geometric scaling factors from the element maps are applied during these projections. The inverse mass matrices are applied only after all volume and surface contributions have been accumulated.

\paragraph{Surface terms} 
Each directional surface contribution is constructed independently by extracting TT traces on the corresponding element faces.
For \(z \in \{x_2,v_{x_1},v_{x_2}\}\), TT traces \(f^-\) and
\(f^+\) are extracted from \(\mathcal{F}^{TT}\) on the two faces normal to
\(z\) by replacing the basis-evaluation matrix in the \(z\)-direction with the
corresponding left or right endpoint evaluation, while the remaining directions
are kept in quadrature form.
For ease of notation, for a given direction \(z\), we write
\(T=T_z^{TT}\), where \(T_z^{TT}\) denotes the TT representation of the
normal transport speed \(T_z\) in the corresponding face-quadrature layout.
The implementation uses the split upwind form
\[
T^\pm = \frac12(T \pm |T|),
\qquad
\widehat{Tf}=T^+f^-+T^-f^+.
\]

Since \(T\) is represented as a TT object, the entrywise nonlinear quantity \(|T|\) cannot be obtained directly using standard TT arithmetic. We therefore approximate \(|T|\) using the TT-cross interpolation procedure described in Section~\ref{subsec:tt-cross}, which constructs a compressed TT approximation from selected entry evaluations without assembling the full tensor. Intermediate TT tensors generated during the flux evaluation are recompressed using TT rounding~\cite{oseledets2011tt} to maintain the prescribed low-rank accuracy tolerance. 
Face fluxes are projected back to modal coefficients by evaluating the basis
functions at the face endpoints in the normal direction and applying the usual
quadrature-weighted test projections in the other two coordinate directions.

\paragraph{Right-hand-side assembly}
After computing the TT representations of the volume and surface contributions
corresponding to $\mathcal{A}_h$, we apply TT rounding and then apply the
one-dimensional inverse mass matrices
$M_{x_2}^{-1}$, $M_{v_{x_1}}^{-1}$, and $M_{v_{x_2}}^{-1}$ to the corresponding
modal cores.
Algorithm~1 summarizes the evaluation of the Vlasov DG right-hand side in TT
format. We denote the resulting right-hand-side operator by
$\operatorname{RHS}^{\mathrm{TT}}$.

\begin{algorithm}[H]
\caption{TT evaluation of the 1D2V Vlasov DG right-hand side}
\KwData{$\mathcal{F}^{\mathrm{TT}}$, field coefficients $\bE_h,\bB_h$, basis and quadrature data, tolerance $\varepsilon_{\mathrm{TT}}$}
\KwResult{$(\partial_t f_h)^{\mathrm{TT}}$}
Compute $\mathcal{F}^{\mathrm{TT,quad}}$ by applying $\Phi_z$ to the modal cores in each phase-space direction\;
Form the three volume residuals by TT products with $v_{x_2}$, $E_{x_1}+v_{x_2}B_{x_3}$, and $E_{x_2}-v_{x_1}B_{x_3}$, followed by weighted test projections\;
\ForEach{$z\in\{x_2,v_{x_1},v_{x_2}\}$}{
  Extract TT traces $f^-$ and $f^+$ on the two faces normal to $z$\;
  Build the normal transport speed $T$ in the face-quadrature TT layout\;
  Approximate $|T|$ by TT-cross interpolation and form $T^+f^-+T^-f^+$ by TT arithmetic\;
  Project the left and right face fluxes back to modal coefficients\;
}
Sum all volume and surface residuals, apply the appropriate sign, and round to $\varepsilon_{\mathrm{TT}}$\;
Apply the one-dimensional inverse mass matrices to the modal cores\;
\end{algorithm}
Algorithm~1 returns the TT representation of $\partial_t f_h$.

\subsection{TT realization of the time update}
The fully discrete TT update follows the full-grid Scheme~\eqref{eq:scheme1},
with the Vlasov right-hand-side evaluations carried out in TT format. In each
time step, the Vlasov DG right-hand side is evaluated twice. The first
evaluation uses the distribution and fields at time $t^n$

\[
\mathcal{R}^{\mathrm{TT},n}
=
\operatorname{RHS}^{\mathrm{TT}}
\left(
\mathcal{F}^{\mathrm{TT},n};
\bE_h^n,\bB_h^n
\right).
\]
The midpoint distribution tensor is then computed as 
\begin{equation}
\mathcal{F}^{\mathrm{TT},n+\frac12}
=
\operatorname{round}\!\left(
\mathcal{F}^{\mathrm{TT},n}
+\frac{\Delta t}{2}\mathcal{R}^{\mathrm{TT},n},
\varepsilon_{\mathrm{TT}}
\right),
\label{eq:tt-midpoint-update}
\end{equation}
where $\mathrm{round}(\cdot,\varepsilon_{\mathrm{TT}})$ denotes TT truncation with tolerance $\varepsilon_{\mathrm{TT}}$.

The midpoint current density is computed by velocity-moment contraction of
\(\mathcal{F}^{\mathrm{TT},n+\frac12}\):
\begin{equation}
\J_h^{n+\frac12}(\bx)
=
\int_{\Omega_v} f_h^{n+\frac12}(\bx,\bv)\,\bv\,d\bv,
\end{equation}
with the velocity quadrature and contractions performed directly in TT form.
This current is used in the Maxwell update for \(\bE_h^{n+1}\) and
\(\bB_h^{n+\frac12}\), as in Scheme~\eqref{eq:scheme1}.

The second Vlasov evaluation uses the midpoint distribution and the time-centered fields,
\[
\mathcal{R}^{\mathrm{TT},n+\frac12}
=
\operatorname{RHS}^{\mathrm{TT}}
\left(
\mathcal{F}^{\mathrm{TT},n+\frac12};
\frac{\bE_h^n+\bE_h^{n+1}}{2},
\bB_h^{n+\frac12}
\right).
\]
The final Vlasov update is then
\begin{equation}
\mathcal{F}^{\mathrm{TT},n+1}
=
\operatorname{round}\!\left(
\mathcal{F}^{\mathrm{TT},n}
+\Delta t\,\mathcal{R}^{\mathrm{TT},n+\frac12},
\varepsilon_{\mathrm{TT}}
\right).
\label{eq:tt-final-update}
\end{equation}
Thus, TT rounding is applied after both the midpoint and final Vlasov updates,
while the electromagnetic fields remain in the full-grid DG representation.

The TT time integration strategy employed here follows the classical
step-and-truncate paradigm commonly used in tensor-train~\cite{kormann2015semilagrangian,rodgers2022adaptive,dolgov_fast-parabolic_2012} and low-rank tensor
methods~\cite{kieri_projection-dynamical_2019,nakao2025reduced} for kinetic plasmas. After each explicit time-update stage, the resulting TT tensor is
recompressed by TT rounding to the prescribed tolerance
$\varepsilon_{\mathrm{TT}}$. This controls the growth of the TT ranks while
maintaining the prescribed approximation accuracy throughout the computation.

\section{Numerical results} \label{sec:sec_num_exp}

In this section, we present numerical experiments to assess the performance of the tensor-train DG formulation. The TT-DG results are compared against a full-grid modal DG solver across the benchmark problems considered in this work.

Both solvers are implemented in MATLAB, using the same underlying DG discretization and explicit time-integration scheme. 
In the TT solver, tensor compression is applied only to the Vlasov equation, while the Maxwell equations are advanced using the full grid DG representation. We report wall-clock runtime, TT ranks, and compression ratio, defined as the ratio of the number of degrees of freedom in the full-grid DG coefficient tensor to the number of stored parameters in the TT representation. 

The FG solver is a parallelized modal DG implementation in which only the computationally intensive Vlasov equation is parallelized across 44 Intel Xeon cores, while the Maxwell solver remains one-dimensional and is not parallelized. The TT implementation is built on the TT-Toolbox \cite{oseledets_tttoolbox}, with TT-cross operations accelerated by the THOR library \cite{adak2026thor}. 

The initial tensor-train representation of the distribution function was constructed using the AMEn-cross routine from the TT-Toolbox. Default parameters were \texttt{nswp}=20, \texttt{zrank}=5, \texttt{zrank2}=3, \texttt{kickrank}=2, \texttt{max\_err\_jumps}=2, \texttt{trunc\_method}=\texttt{'cross'}, and \texttt{tol\_exit}=\(\varepsilon_{\mathrm{TT}}\); the default random rank-2 TT initializer was used. 
The entrywise absolute value \( |T| \), required for the upwind flux splitting, was approximated using the THOR greedy DMRG-cross routine. The computation used the default \texttt{dmrgg\_cross} parameters \texttt{tol} \(= \varepsilon_{\mathrm{TT}}\), \texttt{pivoting} \(=3\), \texttt{vectorized} \(=\texttt{true}\), and \texttt{batch\_size} \(=32\), while \texttt{maxrank} was left unspecified. Internally, the native Fortran implementation initializes the iteration from the default rank-1 all-ones TT tensor.

\subsection{TT layout and phase-ordering study}
\label{subsec:tt_layout_study}

Using the TT layouts introduced in
Section~\ref{subsec:tt-representation-ordering},
we evaluated the impact of tensor layout and phase-space permutation on runtime and compressibility for three $1\text{D}2\text{V}$ benchmark problems considered in this work.

For each problem, we evaluated both the interleaved and grouped layouts together with all six permutations
$
\pi=(\pi_1,\pi_2,\pi_3)
$
of $(x_2,v_{x_1},v_{x_2})$. Layouts were ranked primarily by mean wall-clock time, with run-averaged TT compression used as a secondary metric. All three layout-study campaigns were carried out to the same final time, $T=1$.

\begin{table}[H]
\centering
\small
\begin{tabular}{llll}
\toprule
Problem & Tested grids & Best layout & Comparison with worst tested layout \\
\midrule
Streaming Weibel instability &
\makecell[l]{$(30,30,30)$\\$(60,60,60)$\\$(120,120,120)$} &
\makecell[l]{interleaved\\$\pi=(v_{x_2},v_{x_1},x_2)$} &
\makecell[l]{worst: grouped, $\pi=(x_2,v_{x_1},v_{x_2})$\\$1.99\times$--$4.73\times$ faster\\$2.58\times$--$3.47\times$ higher mean compression} \\
Weak Landau damping &
\makecell[l]{$(32,32,32)$\\$(32,64,64)$\\$(32,128,128)$} &
\makecell[l]{interleaved\\$\pi=(v_{x_2},x_2,v_{x_1})$} &
\makecell[l]{worst: grouped, $\pi=(x_2,v_{x_1},v_{x_2})$\\$1.91\times$--$1.96\times$ faster\\$6.14\times$--$8.58\times$ higher mean compression} \\
Two-stream instability &
\makecell[l]{$(32,32,32)$\\$(32,64,64)$\\$(32,128,128)$} &
\makecell[l]{interleaved\\$\pi=(v_{x_1},x_2,v_{x_2})$} &
\makecell[l]{worst: grouped, $\pi=(x_2,v_{x_1},v_{x_2})$\\$2.00\times$--$16.45\times$ faster\\$4.49\times$--$12.28\times$ higher mean compression} \\
\bottomrule
\end{tabular}
\caption{Summary of TT layout sweeps for the $1\text{D}2\text{V}$ benchmark problems, all performed to final time $T=1$. The final column reports the speedup and mean compression improvement of the best configuration relative to the slowest tested layout across the listed grids. All TT layouts exhibit compression and perform faster than their full-grid counterparts.}
\label{tab:tt_layout_recommendations}
\end{table}

Table~\ref{tab:tt_layout_recommendations} summarizes the best performing configurations for each benchmark. Across all cases, the interleaved layout consistently outperformed the grouped layout, yielding both lower runtime and higher average compression. The performance gap generally increased with velocity resolution, indicating that tensor ordering becomes increasingly important as the dimensional cost of the DG discretization grows. Even the least favorable TT layouts outperformed the corresponding full-grid solver on the tested problems, indicating that the benefits of TT compression are robust to the choice of ordering.

The optimal permutation $\pi$ was problem dependent. For the streaming Weibel instability, the most consistent performance across grids was obtained with 
$
\pi=(v_{x_2},v_{x_1},x_2).
$
Weak Landau damping favored
$
\pi=(v_{x_2},x_2,v_{x_1}),
$
while the best performance for two-stream instability was obtained with
$
\pi=(v_{x_1},x_2,v_{x_2}).
$
The improved performance of the interleaved layout is likely due to the placement of element and modal indices associated with the same coordinate direction in adjacent TT cores, which reduces intermediate TT-rank growth during operator application.
These configurations are used in the experiments reported below.
\subsection{Streaming Weibel instability}

We consider the streaming Weibel instability benchmark for the Vlasov--Maxwell system as presented in Cheng et al. \cite{cheng2014energy}. In addition to providing a representative coupled kinetic–electromagnetic instability problem, this benchmark enables a stringent time-reversibility test for assessing the accuracy of the TT-DG formulation. The spatial domain is periodic with length $\Omega = 2\pi/k_0$, where $k_0 = 0.2$, and the velocity domain is truncated to $v_{x_1},v_{x_2}\in[-1.2,1.2]$. The initial distribution is
\[
f_0(x_2,v_{x_1},v_{x_2})
=
\left[
(1-\delta)\exp\!\left(-\frac{(v_{x_1}+v_{x_2,0})^2}{\beta}\right)
+
\delta\exp\!\left(-\frac{(v_{x_1}-v_{x_1,0})^2}{\beta}\right)
\right]
\frac{\exp\!\left(-v_{x_2}^2/\beta\right)}{\pi\beta},
\]
with $\beta=0.01$, $\delta=0.5$, and $v_{x_1,0}=v_{x_2,0}=0.3$. The electromagnetic fields are initialized by
\[
E_1(x_2,0)=E_2(x_2,0)=0,
\qquad
B_3(x_2,0)=b\sin(k_0x_2),
\]
with $b=10^{-3}$.

We use a modal DG discretization of degree $p=1$ on uniform grids with
\[
N_{x_2}=N_{v_{x_1}}=N_{v_{x_2}}
\in
\{30,40,50,60,80,120,160\}.
\]
Time integration uses
\[
\Delta t
=
\mathrm{CFL}\,
\min(\Delta x_2,\Delta v_{x_1},\Delta v_{x_2}),
\qquad
\mathrm{CFL}=0.15,
\]
with final time $T=1$. The TT solver uses truncation tolerance $\varepsilon_{\mathrm{TT}}=10^{-6}$ and the interleaved layout with
$
\pi=(v_{x_2},v_{x_1},x_2)
$
identified in Section~\ref{subsec:tt_layout_study}.

To assess accuracy, we employ the time-reversal test of Cheng et al. \cite{cheng2014energy}. The system is evolved from $t=0$ to $t=T$, after which the velocity is reversed, $\bv\mapsto-\bv$, and the magnetic field sign is negated before evolving for an additional interval of length $T$. Since the Vlasov--Maxwell system is time reversible, the exact solution at time $2T$ coincides with the initial condition.

Let $f_h^{\mathrm{rev}}$ denote the numerical solution at time $2T$. We define the reversal error by
\[
e_{\mathrm{rev}}^f
=
\left\|
f_h^{\mathrm{rev}}(\cdot,\cdot,\cdot,2T)
-
f_h(\cdot,\cdot,\cdot,0)
\right\|_{L^2(\Omega)}.
\]
The DG solution is reconstructed at tensor-product quadrature points and evaluated against the prescribed initial distribution $f_0$ using elementwise Gauss quadrature of order $2p+1$.

\begin{table}[H]
\centering
\small
\begin{tabular}{ccccc}
\toprule
$N$ & FG error in $f$ & FG rate & TT error in $f$ & TT rate \\
\midrule
30 & $1.05\mathrm{E}{-02}$ & --    & $1.00\mathrm{E}{-02}$ & -- \\
40 & $7.21\mathrm{E}{-03}$ & $1.30$ & $6.92\mathrm{E}{-03}$ & $1.28$ \\
50 & $5.02\mathrm{E}{-03}$ & $1.62$ & $4.78\mathrm{E}{-03}$ & $1.66$ \\
60 & $3.65\mathrm{E}{-03}$ & $1.75$ & $3.45\mathrm{E}{-03}$ & $1.80$ \\
80 & $2.15\mathrm{E}{-03}$ & $1.84$ & $1.84\mathrm{E}{-03}$ & $2.18$ \\
120 & $9.94\mathrm{E}{-04}$ & $1.90$ & $8.52\mathrm{E}{-04}$ & $1.90$ \\
160 & $5.89\mathrm{E}{-04}$ & $1.82$ & $4.89\mathrm{E}{-04}$ & $1.93$ \\
\bottomrule
\end{tabular}
\caption{Time-reversal error and observed convergence rate for the distribution function $f$ in the full-grid and TT solvers.}
\label{tab:time_reversal_f_errors}
\end{table}

Table~\ref{tab:time_reversal_f_errors} reports the reversal error and observed convergence rates. Both the FG and TT solvers exhibit approximately second-order convergence, with comparable reversal errors across all reported grids. Exact reversibility is not expected due to the time discretization and, in the TT formulation, additional truncation and cross-interpolation errors. Despite these effects, both methods maintain good long-time reversibility behavior. 

\begin{figure}[H]
\centering
\subfigure[Measured runtime versus uniform phase-space resolution for the full-grid and TT solvers.]{
\includegraphics[width=0.48\textwidth]{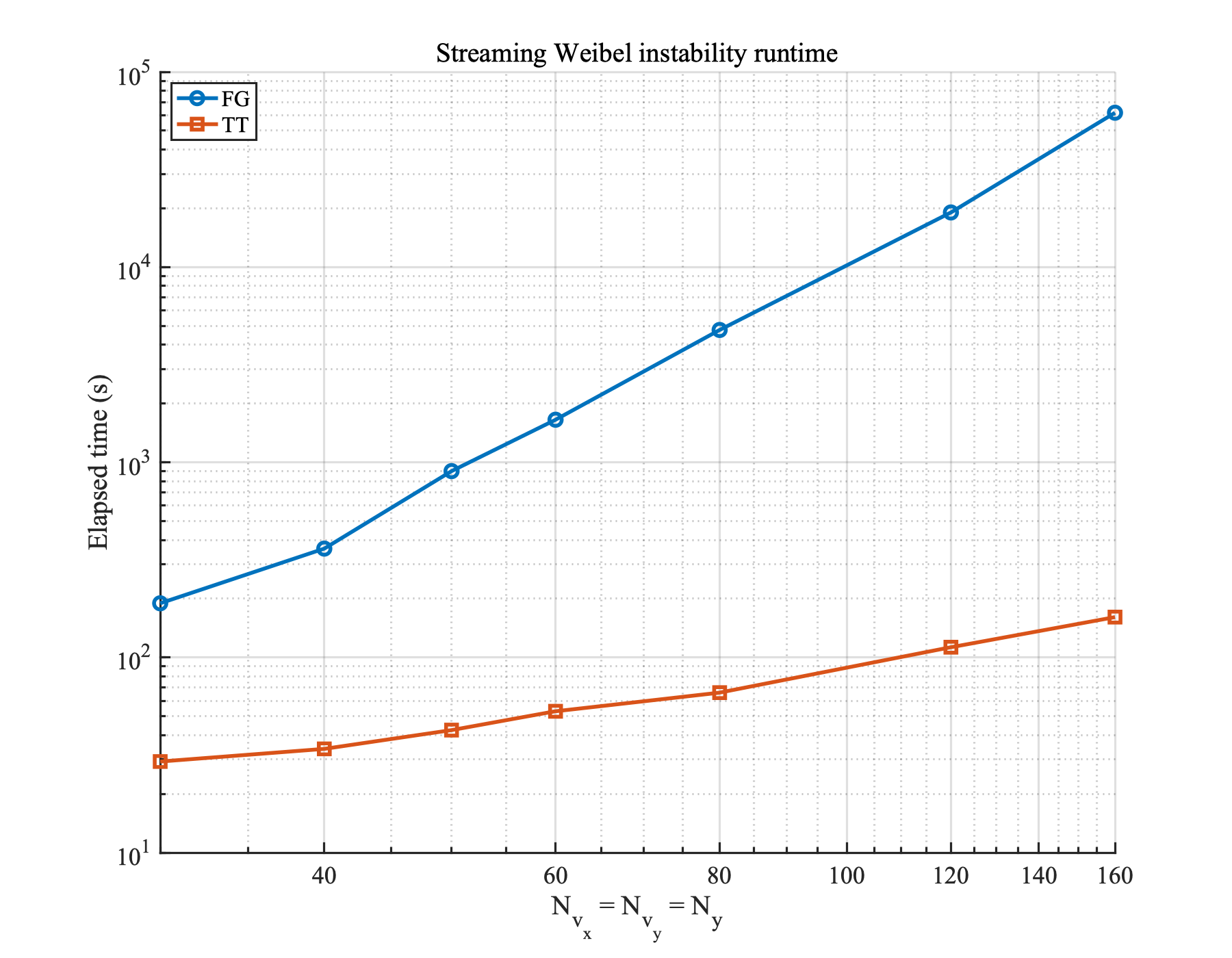}
\label{fig:time_reversal_runtime}
}
\hfill
\subfigure[TT compression ratio versus uniform phase-space resolution for the time-reversal benchmark.]{
\includegraphics[width=0.48\textwidth]{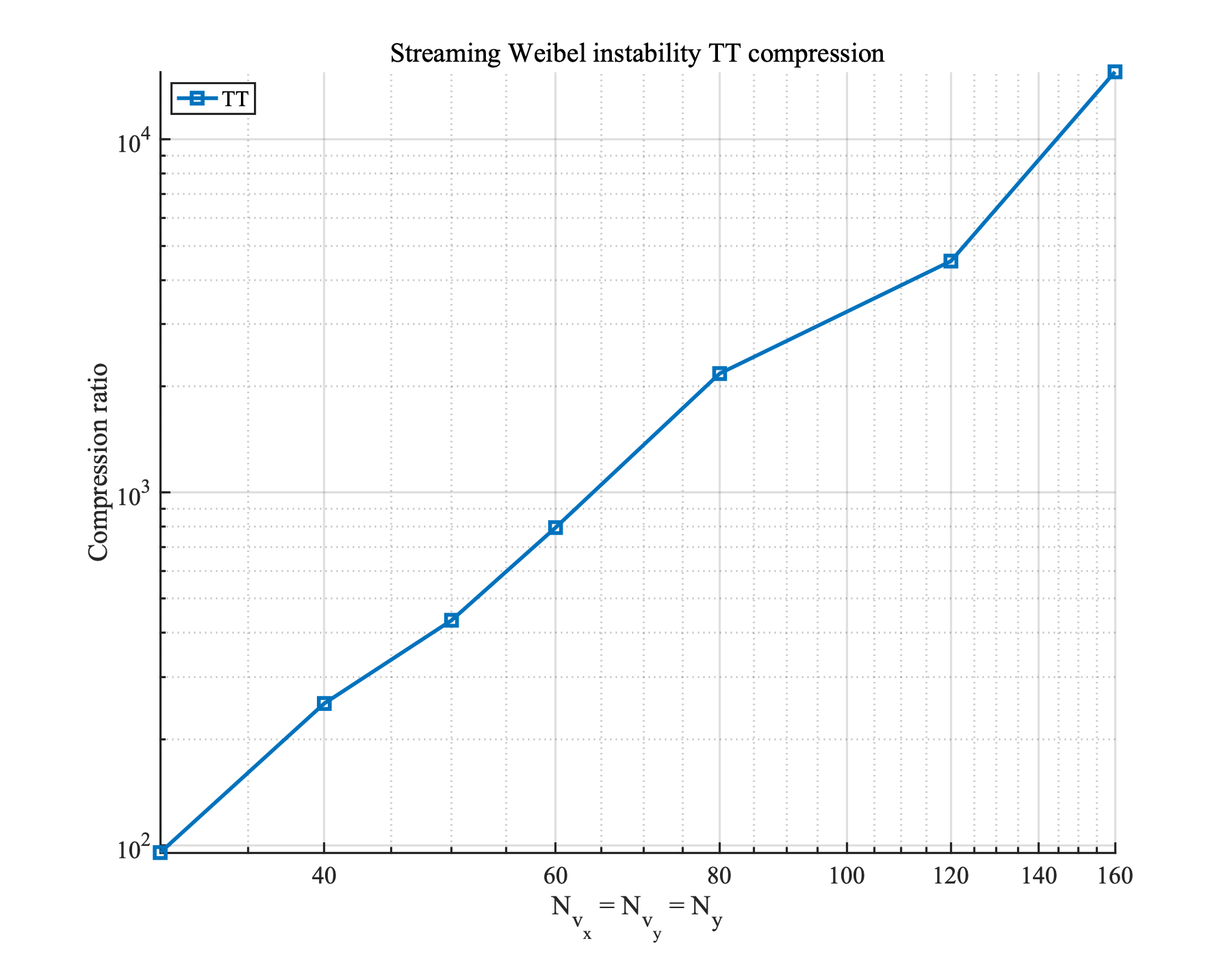}
\label{fig:time_reversal_compression}
}
\caption{Runtime and TT compression at $T=1$ for the streaming Weibel time-reversal benchmark. The TT solver exhibits increasing speedup and compression with phase-space refinement.}
\label{fig:time_reversal_performance}
\end{figure}

Figure~\ref{fig:time_reversal_performance} summarizes the runtime and compression behavior. The TT solver is consistently faster than the FG method, with speedups increasing from $6.5\times$ at $N=30$ to $3.84\times10^2$ at $N=160$. At the same time, the TT compression ratio grows monotonically from $9.53\times10^1$ to $1.56\times10^4$, indicating that the streaming Weibel benchmark remains highly compressible as the full-grid DG cost increases under refinement.
\subsection{Weak Landau damping}

We consider the weak Landau damping benchmark as studied in \cite{kormann2015semilagrangian} for the Vlasov–Maxwell system with initial data on the electrostatic invariant submanifold. This problem provides a standard test of linear kinetic damping and long-time accuracy in collisionless plasma simulations. The spatial domain is periodic with $x_2\in[0,4\pi]$, corresponding to the fundamental wavenumber $k_0=0.5$, while the velocity domain is truncated to $v_{x_1},v_{x_2}\in[-6,6]$. The initial distribution is the perturbed Maxwellian
\[
f_0(x_2,v_{x_1},v_{x_2})
=
\frac{1}{2\pi}
\exp\!\left(
-\frac{v_{x_1}^2+v_{x_2}^2}{2}
\right)
\left(1+\alpha\cos(k_0x_2)\right),
\]
with perturbation amplitude $\alpha=0.01$. The electromagnetic fields are initialized by
\[
E_1(x_2,0)=0,
\qquad
E_2(x_2,0)=\frac{\alpha}{k_0}\sin(k_0x_2),
\qquad
B_3(x_2,0)=0.
\]

The phase space is discretized using a modal DG method with polynomial degree $p=1$ on uniform grids with
\[
(N_{x_2},N_{v_{x_1}},N_{v_{x_2}})
\in
\{
(32,32,32),
(32,48,48),
(32,64,64),
(32,80,80),
(32,128,128),
(32,256,256)
\}.
\]
Time integration uses the CFL timestep
\[
\Delta t
=
\mathrm{CFL}\,
\min(\Delta v_{x_1},\Delta v_{x_2}),
\qquad
\mathrm{CFL}=0.05,
\]
with final time $T=20$. The TT solver uses truncation tolerance $\varepsilon_{\mathrm{TT}}=10^{-6}$ together with the interleaved layout and permutation
$
\pi=(v_{x_2},x_2,v_{x_1})
$
identified in Section~\ref{subsec:tt_layout_study}.

To quantify damping, we track the fundamental electric-field Fourier mode
\[
\hat{E}_{k_0}(t)
=
\frac{1}{L_{x_2}}
\int_0^{L_{x_2}}
E_2(x_2,t)e^{-ik_0x_2}\,dx_2.
\]
The theoretical damping rate for this problem is $-.153359$ \cite{ARBER2002339}.

\begin{figure}[t]
\centering
\includegraphics[width=.88\textwidth]{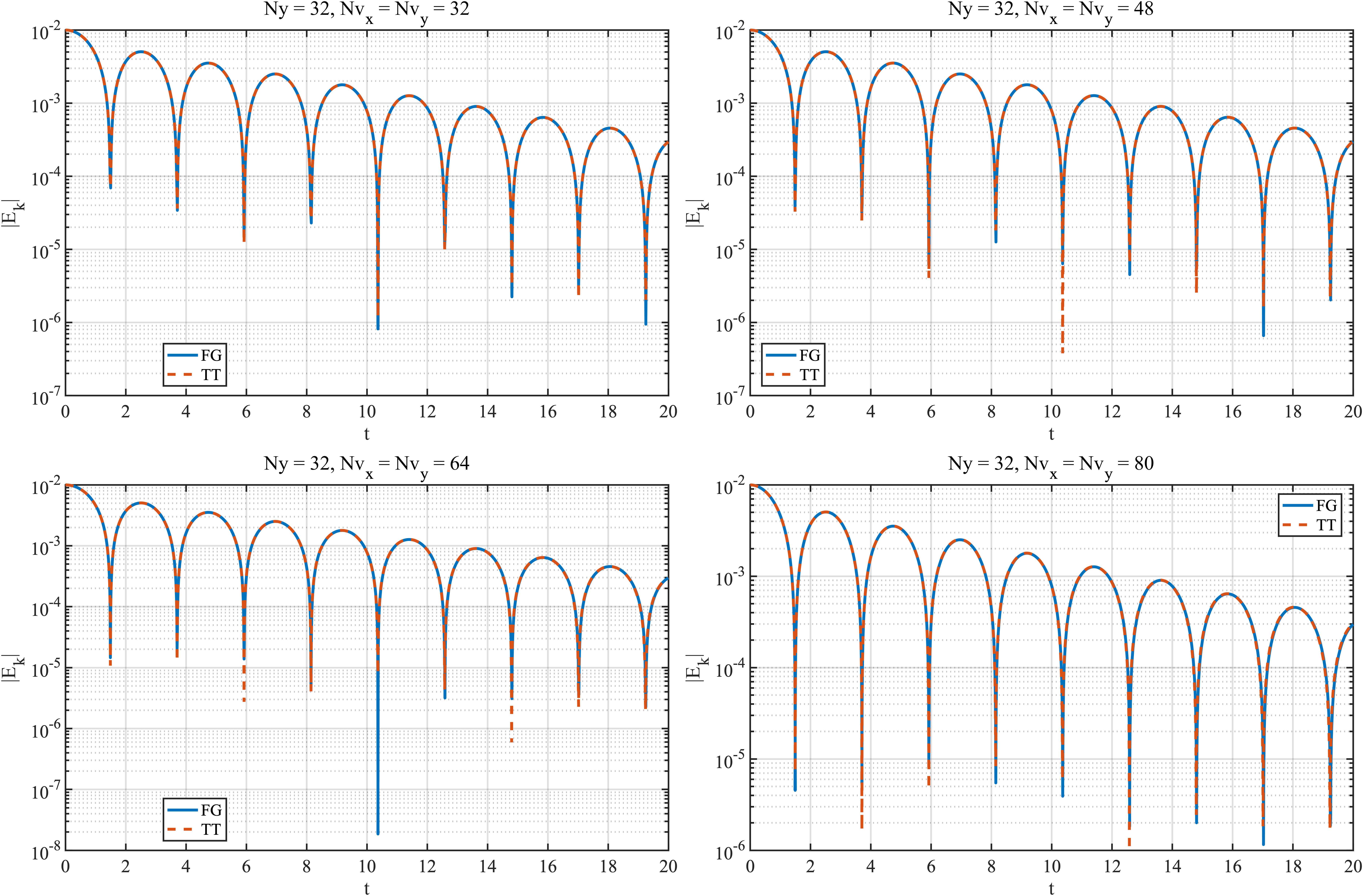}
\caption{Evolution of the fundamental electric-field Fourier mode for weak Landau damping, comparing FG and TT solutions across the shared completed resolutions.}
\label{fig:landau_ek}
\end{figure}

Figure~\ref{fig:landau_ek} compares the evolution of $\hat{E}_{k_0}(t)$ for the FG and TT formulations. Over the shared completed resolutions
$
N_{v_{x_1}}=N_{v_{x_2}}\in\{32,48,64,80\},
$
the curves are visually indistinguishable throughout the linear damping interval, indicating that the TT truncation preserves the damping dynamics accurately.

\begin{table}[H]
\centering
\small
\begin{tabular}{ccc}
\toprule
$N_{v_{x_1}}=N_{v_{x_2}}$ & FG damping rate & TT damping rate \\
\midrule
32  & $-0.15468$ & $-0.15475$ \\
48  & $-0.15457$ & $-0.15435$ \\
64  & $-0.15452$ & $-0.15448$ \\
80  & $-0.15449$ & $-0.15436$ \\
128 & --         & $-0.15374$ \\
256 & --         & $-0.15771$ \\
\bottomrule
\end{tabular}
\caption{Fitted damping rates for weak Landau damping. The FG study completed through $N_{v_{x_1}}=N_{v_{x_2}}=80$, whereas the TT study continued through $N_{v_{x_1}}=N_{v_{x_2}}=256$.}
\label{tab:landau_damping_rates}
\end{table}

Table~\ref{tab:landau_damping_rates} lists the fitted damping rates. For the shared completed resolutions, the FG and TT rates agree to at least three digits and remain close to the theoretical value. The FG runs became computationally impractical beyond
$
N_{v_{x_1}}=N_{v_{x_2}}=80,
$
whereas the TT formulation completed simulations through $256^2$ velocity resolution. At
$
N_{v_{x_1}}=N_{v_{x_2}}=256,
$
the TT damping rate shows a modest deterioration, likely due to accumulated tensor truncation effects at very high compression ratios.

\begin{figure}[H]
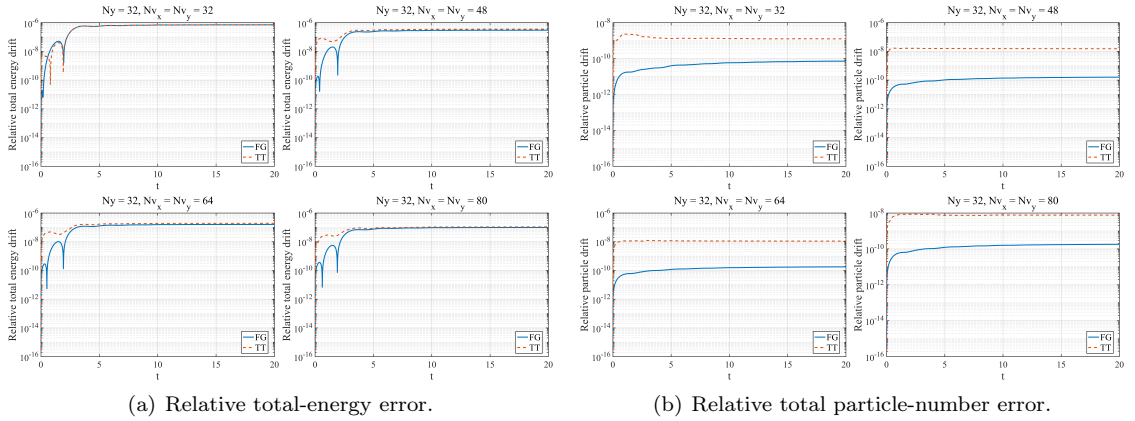

\centering
\subfigure[Relative total-energy error.]{
\includegraphics[width=0.48\textwidth]{landau_fg_vs_tt_total_energy_drift.png}
}
\hfill
\subfigure[Relative total particle-number error.]{
\includegraphics[width=0.48\textwidth]{landau_fg_vs_tt_relative_particle_drift.png}
}
\caption{Relative total-energy and particle-number errors for weak Landau damping}
\label{fig:landau_conservation}
\end{figure}

Figure~\ref{fig:landau_conservation} shows the relative conservation diagnostics. Both formulations maintain small total-energy and particle-number errors over the simulation interval, indicating that the TT formulation preserves the dominant conservation properties of the underlying DG discretization despite the additional truncation and cross-interpolation steps. 

\begin{figure}[H]
\centering
\includegraphics[width=0.88\textwidth]{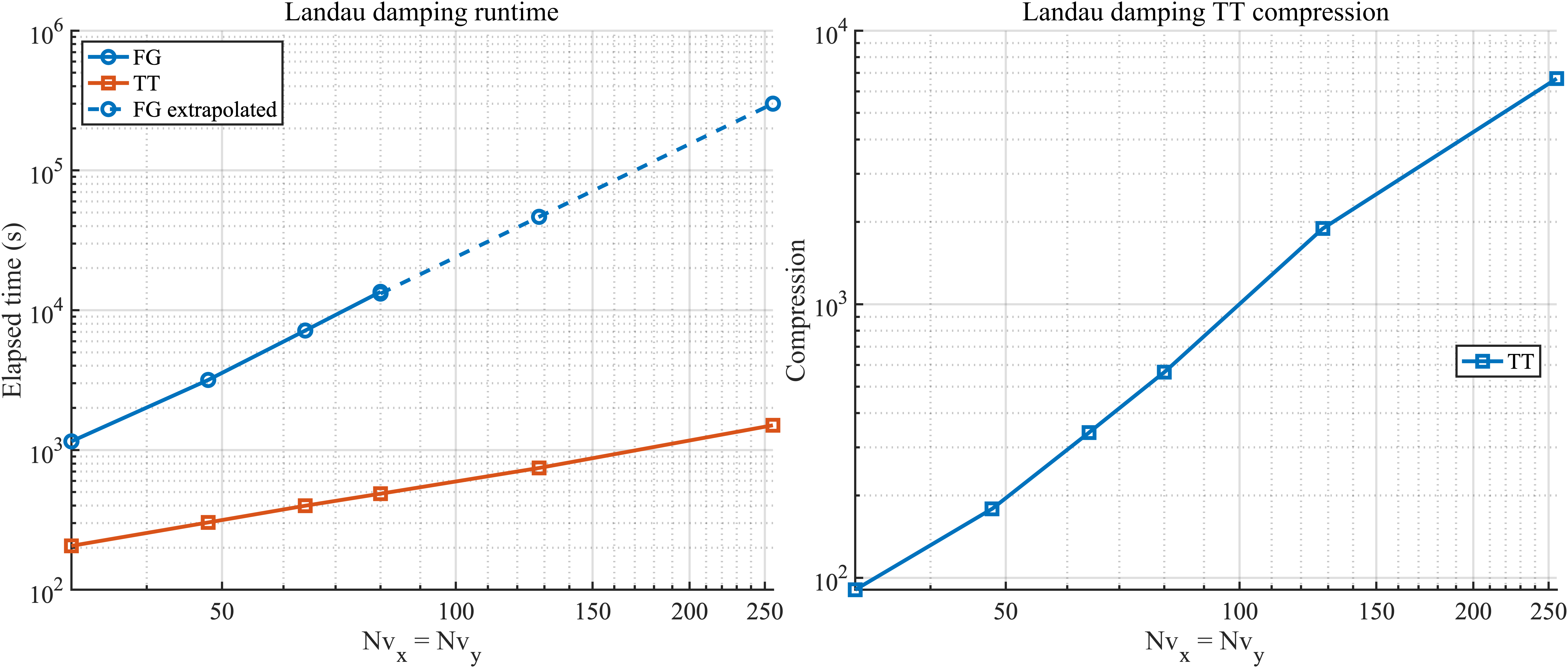}
\caption{Runtime and TT compression at $T=20$ versus velocity resolution for weak Landau damping.}
\label{fig:landau_performance}
\end{figure}

Figure~\ref{fig:landau_performance} summarizes the computational performance. The TT solver substantially reduces runtime relative to the FG formulation while enabling simulations at resolutions that were not computationally tractable in the full-grid setting. At the same time, the TT compression ratio increases strongly with refinement indicating that the weak Landau damping solution remains highly compressible throughout the linear regime.

\begin{figure}[H]
\centering
\includegraphics[width=0.78\textwidth]{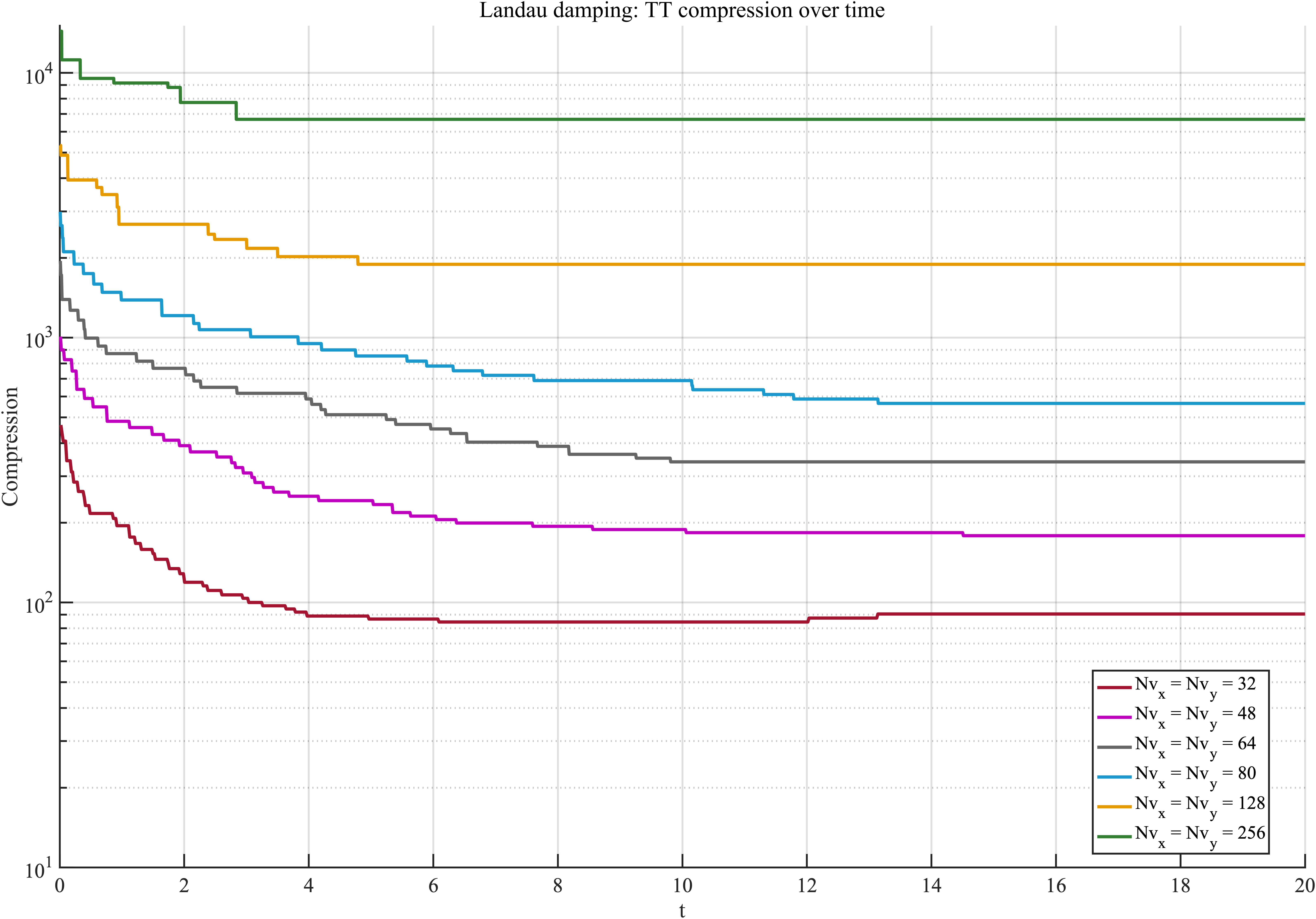}
\caption{TT compression ratio over time for weak Landau damping.}
\label{fig:landau_tt_compression_time}
\end{figure}

Figure~\ref{fig:landau_tt_compression_time} shows the temporal evolution of the TT compression ratio. After a brief initial transient, the compression stabilizes into a long plateau for all reported resolutions, demonstrating sustained compressibility over the full damping interval.

\subsection{Two-stream instability}
We consider the magnetically induced two-stream instability benchmark for the Vlasov–Maxwell system studied in \cite{crouseilles2015hamiltonian}. This problem provides a stringent test of nonlinear phase-space filamentation and the ability of the TT formulation to remain effective in the presence of rapidly developing fine-scale structure.

The initial distribution consists of two counter-propagating beams in the $v_{x_2}$ direction,
\[
f_0(x_2,v_{x_1},v_{x_2})
=
\frac{1}{2\pi\beta}
\exp\!\left(-\frac{v_{x_1}^2}{\beta}\right)
\left[
\exp\!\left(-\frac{(v_{x_2}-v_{x_2,0})^2}{\beta}\right)
+
\exp\!\left(-\frac{(v_{x_2}+v_{x_2,0})^2}{\beta}\right)
\right],
\]
with $\beta=2\times10^{-3}$ and $v_{x_2,0}=0.2$. The electromagnetic fields are initialized by
\[
E_1(x_2,0)=E_2(x_2,0)=0,
\qquad
B_3(x_2,0)=\alpha\sin(x_2),
\]
with perturbation amplitude $\alpha=10^{-3}$.

The phase space is discretized using a modal DG method with polynomial degree $p=1$ on uniform grids with
\[
(N_{x_2},N_{v_{x_1}},N_{v_{x_2}})
\in
\{
(32,32,32),
(32,48,48),
(32,64,64),
(32,80,80),
(32,128,128),
(32,256,256)
\}.
\]
Simulations are run to final time $T=60$. The TT solver uses truncation tolerance $\varepsilon_{\mathrm{TT}}=10^{-8}$ together with the interleaved layout and permutation
$
\pi=(v_{x_1},x_2,v_{x_2})
$
identified in Section~\ref{subsec:tt_layout_study}.

The two-stream instability is substantially more challenging than the weak Landau damping problem due to its strongly nonlinear evolution, rapid phase-space filamentation, and the formation of fine-scale vortex structures after saturation. These effects reduce compressibility and make the problem particularly demanding for tensor-based approximations. In practice, larger truncation tolerances produced qualitatively different dynamics relative to the full-grid reference solution, causing premature nonlinear growth and vortex structures with incorrect spatial scales. The tighter tolerance $\varepsilon_{\mathrm{TT}} = 10^{-8}$ was required to recover physically consistent behavior.

\begin{figure}[H]
\centering
\includegraphics[width=0.88\textwidth]{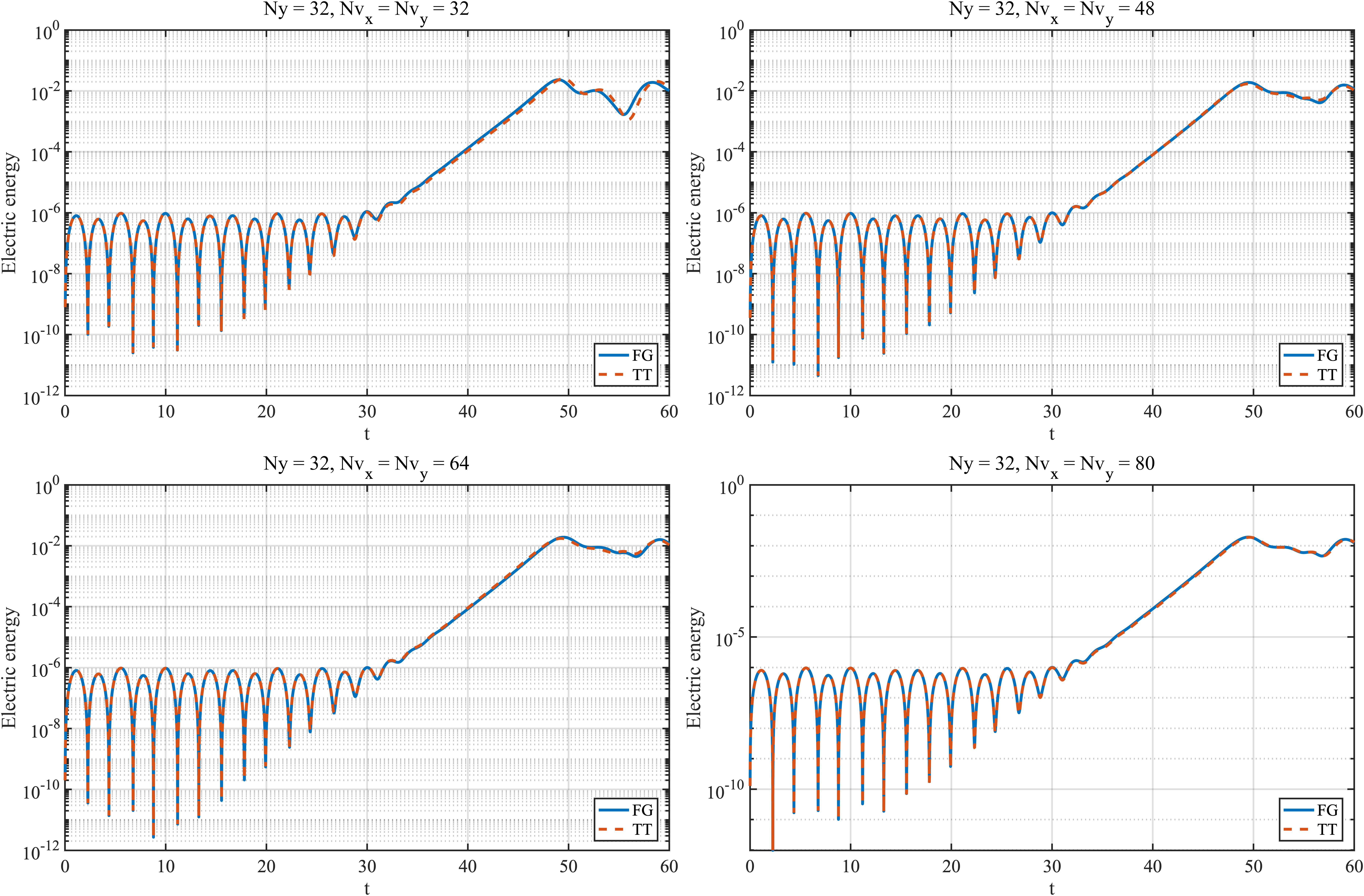}
\caption{Electric energy for the two-stream instability, comparing FG and TT simulations across the shared resolutions.}
\label{fig:two_stream_electric_energy}
\end{figure}

Figure~\ref{fig:two_stream_electric_energy} compares the electric-energy histories for the FG and TT simulations. Over the shared completed resolutions
$
N_{v_{x_1}}=N_{v_{x_2}}\in\{32,48,64,80\},
$
the TT curves closely track the FG results through the linear growth phase, saturation, and subsequent nonlinear evolution, indicating that the TT representation accurately captures the dominant instability dynamics despite the challenging nonlinear behavior.

\begin{figure}[H]
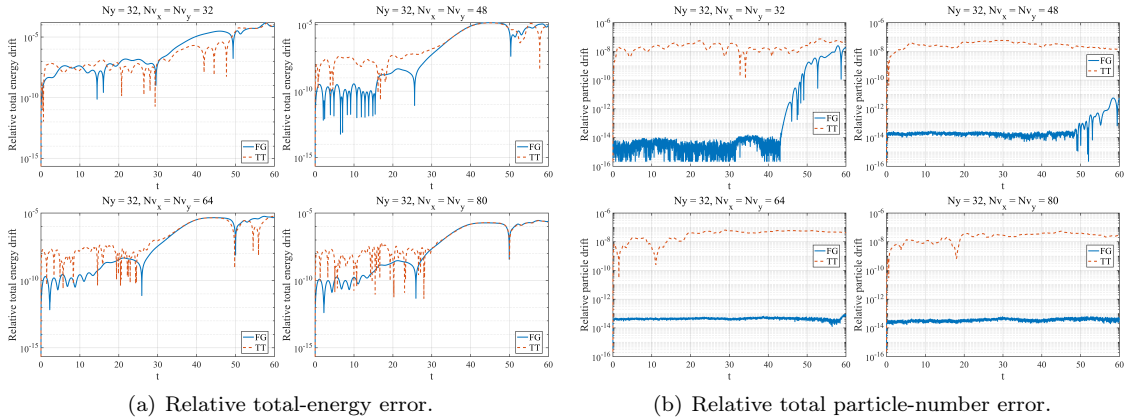

\centering
\subfigure[Relative total-energy error.]{
\includegraphics[width=0.48\textwidth]{two_stream_fg_vs_tt_total_energy_drift.png}
}
\hfill
\subfigure[Relative total particle-number error.]{
\includegraphics[width=0.48\textwidth]{two_stream_fg_vs_tt_relative_particle_drift.png}
}
\caption{Relative conservation diagnostics for the two-stream instability.}
\label{fig:two_stream_conservation}
\end{figure}

Figure~\ref{fig:two_stream_conservation} shows the relative conservation diagnostics. Both formulations maintain small total-energy and particle-number errors throughout the simulation interval, with the errors decreasing under refinement on the shared completed grids. This indicates that the TT compression preserves the dominant conservation properties of the underlying DG discretization even in the strongly nonlinear regime.

\begin{figure}[H]
\centering
\subfigure[Full-grid, $(32,80,80)$.]{
\includegraphics[width=0.48\textwidth]{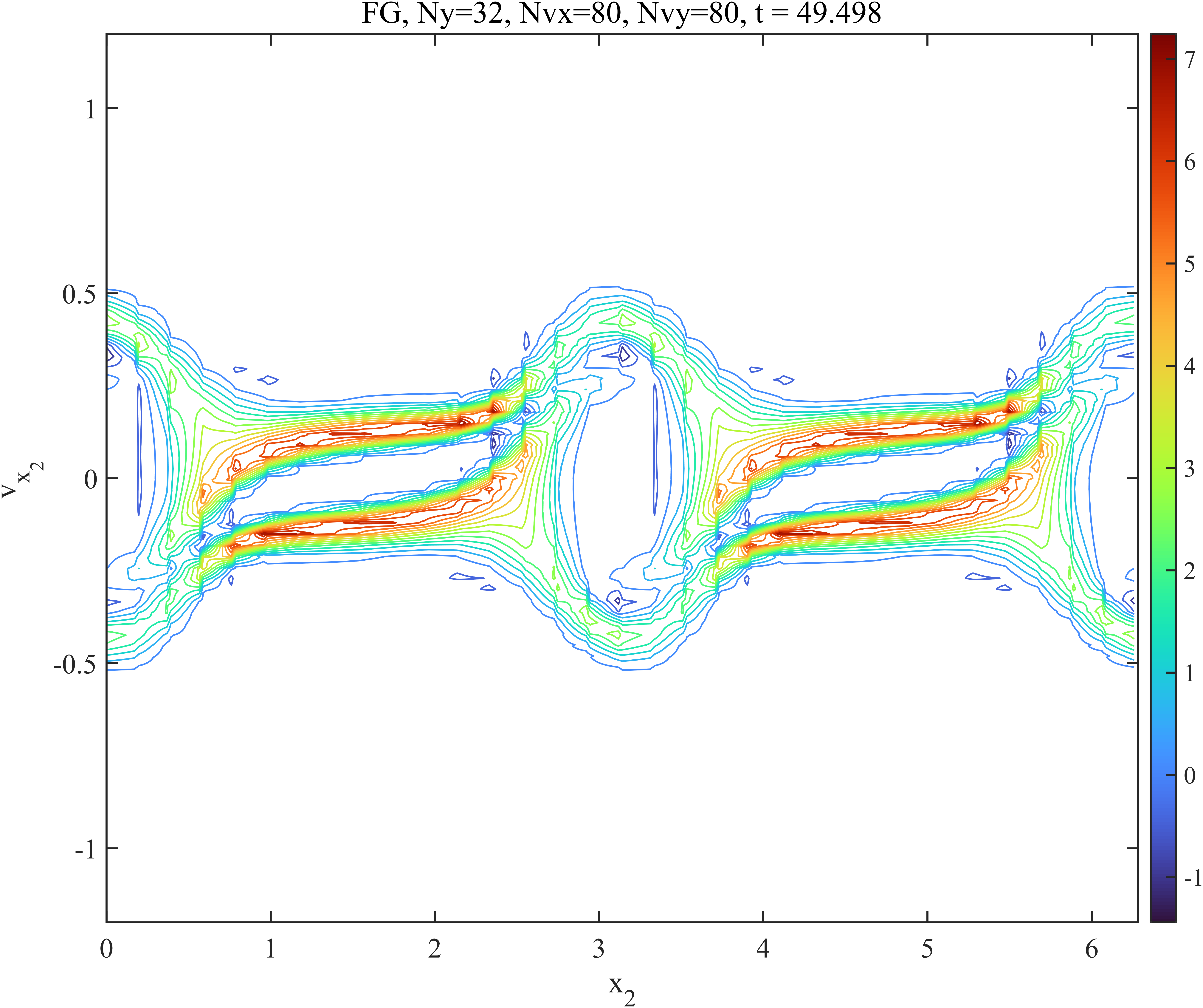}
}
\hfill
\subfigure[TT, $(32,80,80)$.]{
\includegraphics[width=0.48\textwidth]{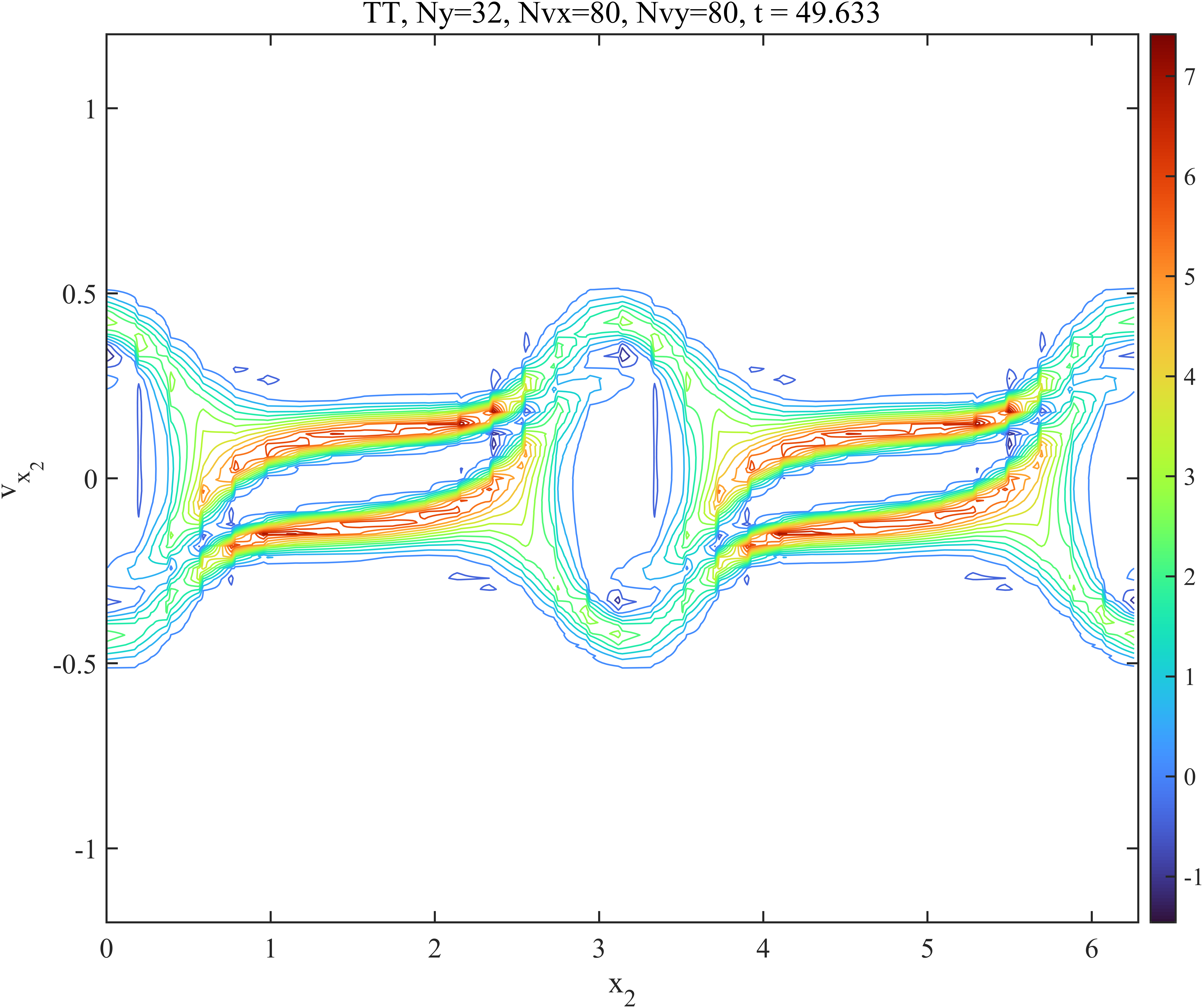}
}
\caption{ Phase-space contour projections in the $(x_2,v_{x2})$ plane near peak electric response for the two-stream instability at resolution $(32,80,80)$.}
\label{fig:two_stream_peak_snapshot}
\end{figure}

To compare the phase-space structure directly, Figure~\ref{fig:two_stream_peak_snapshot} shows $(x_2,v_{x_2})$ contour projections near the time of peak electric response for the representative grid $(32,80,80)$. The TT solution reproduces the filamentation patterns and beam deformation observed in the FG solution, providing qualitative confirmation that the compressed representation remains faithful during the most dynamically active stage of the evolution.

\begin{figure}[H]
\centering
\includegraphics[width=0.88\textwidth]{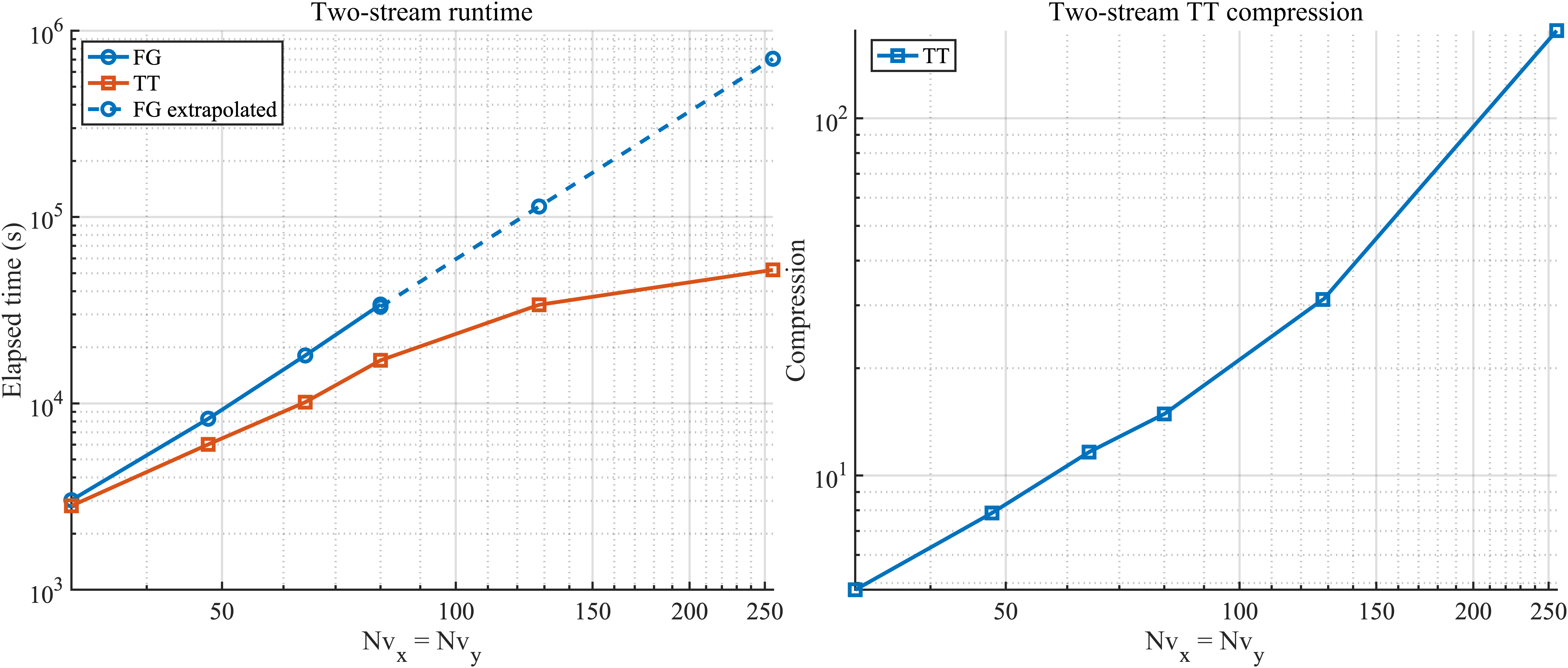}
\caption{Runtime and TT compression at $T=60$ versus velocity resolution for the two-stream instability.}
\label{fig:two_stream_performance}
\end{figure}

Figure~\ref{fig:two_stream_performance} summarizes the computational performance. Compared with the weak Landau damping problem, the attainable TT compression is substantially lower due to the progressive generation of fine-scale nonlinear structure. Nevertheless, the TT formulation still reduces runtime relative to the FG solver on the shared completed grids while enabling simulations at finer resolutions that were not computationally tractable in the full-grid setting. The compression ratio continues to increase with refinement, reaching its largest values on the $128^2$ and $256^2$ velocity grids.

\begin{figure}[H]
\centering
\includegraphics[width=0.78\textwidth]{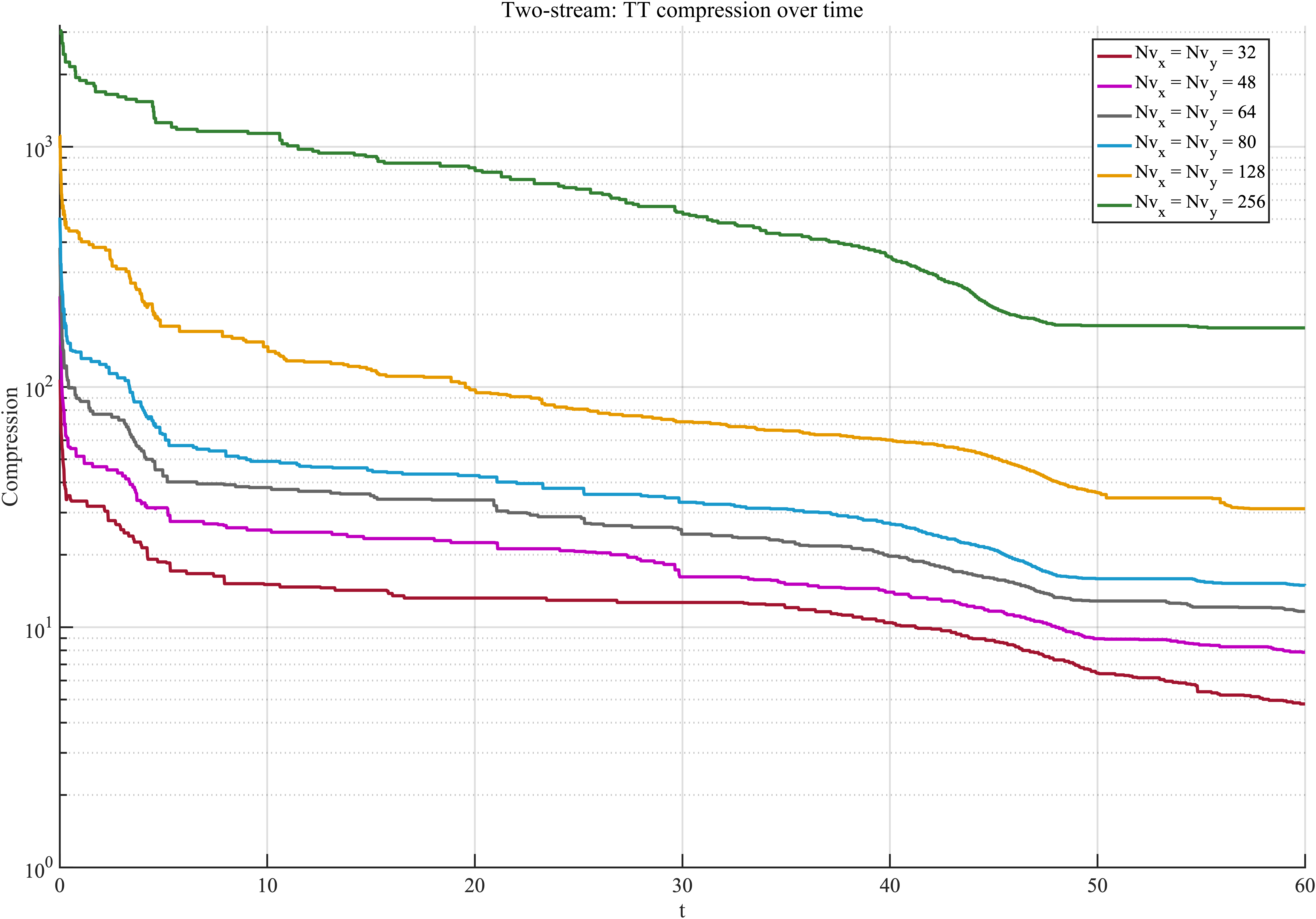}
\caption{Temporal evolution of the TT compression ratio for the two-stream instability.}
\label{fig:two_stream_tt_compression_time}
\end{figure}

Figure~\ref{fig:two_stream_tt_compression_time} shows the temporal evolution of the TT compression ratio. Unlike the weak Landau damping case, the compression decreases gradually over time as nonlinear filamentation generates increasingly fine phase-space structure. Even so, the solution remains compressed throughout the simulation interval, demonstrating that the TT representation remains effective even for this strongly nonlinear benchmark.

\section{Conclusion} \label{sec:sec_conclusion}

We have presented a tensor-train discontinuous Galerkin (TT-DG) formulation for the Vlasov--Maxwell system that combines a modal DG discretization with low-rank tensor representations of the phase-space solution and discrete operators. The formulation exploits the tensor-product structure of the DG discretization to perform quadrature, differentiation, flux evaluation, and time integration directly in compressed form without explicitly constructing full multidimensional arrays.

The proposed method was evaluated on several standard 1D2V Vlasov--Maxwell benchmarks, including streaming Weibel instability, weak Landau damping, and two-stream instability. Across these problems, the TT formulation reproduced the accuracy and conservation behavior of the underlying full-grid DG discretization while substantially reducing memory usage and runtime. For the streaming Weibel and weak Landau damping problems, the TT representation remained highly compressible under refinement, producing compression ratios exceeding $10^3$ together with significant speedups relative to the full-grid solver. For the strongly nonlinear two-stream instability, compressibility decreased as fine-scale filamentation developed, but the TT formulation nevertheless remained effective and accurately captured the dominant instability dynamics and phase-space structure.

We also investigated the influence of tensor layout and phase-space ordering on performance. Across all benchmark problems considered, interleaved tensor layouts consistently outperformed grouped layouts, yielding both lower runtime and improved compression behavior.
These results demonstrate that tensor-train representations provide an effective approach for reducing the computational cost of deterministic kinetic plasma simulations while retaining the favorable accuracy and conservation properties of DG discretizations. 

Although the numerical experiments presented here focus on the 1D2V Vlasov--Maxwell system with polynomial degree $p=1$, the proposed TT-DG framework readily accommodates higher-order polynomial approximations through the same tensor-product construction. Existing full-grid DG studies from Cheng et al. \cite{cheng2014energy} have demonstrated optimal convergence behavior for $p=2$ and $p=3$, motivating a future investigation of the interaction between high-order accuracy and TT truncation. Beyond higher-order discretizations, extending the present approach to higher-dimensional phase spaces represent a natural
next step, where the computational advantages of low-rank compression are expected to become increasingly significant.

\section*{Declaration of generative AI and AI-assisted technologies in the manuscript preparation process}

During the preparation of this work, the authors used OpenAI Codex to assist with code prototyping and implementation and used ChatGPT to assist with language editing and manuscript refinement. After using these tools, the authors tested, reviewed and edited the content as needed and take full responsibility for the content of the published article.


\section*{Acknowledgments}
The authors gratefully acknowledge the support of the Laboratory Directed Research and Development (LDRD) program of Los Alamos National Laboratory under projects 20240705ERil, as well as ISTI Rapid Response support under project 20258349CT-IST and Institutional Computing resources. WB was supported by LDRD under project 20251151PRD1. Los Alamos National Laboratory is operated by Triad National Security, LLC, for the National Nuclear Security Administration of the U.S. Department of Energy under Contract No.~89233218CNA000001.

This research was sponsored by the U.S. Department of Energy Office of Science and the U.S. Air Force Office of Scientific Research.  The authors also gratefully acknowledge support from the LDRD program at Sandia National Laboratories, under project 233973.  Sandia National Laboratories is a multi-mission laboratory managed and operated by National Technology \& Engineering Solutions of Sandia, LLC (NTESS), a wholly owned subsidiary of Honeywell International Inc., for the U.S. Department of Energy’s National Nuclear Security Administration (DOE/NNSA) under contract DE-NA0003525. This written work is coauthored by an employee of NTESS. The employee, not NTESS, owns the right, title and interest in and to the written work and is responsible for its contents. Any subjective views or opinions that might be expressed in the written work do not necessarily represent the views of the U.S. Government. The publisher acknowledges that the U.S. Government retains a non-exclusive, paid-up, irrevocable, world-wide license to publish or reproduce the published form of this written work or allow others to do so, for U.S. Government purposes. The DOE will provide public access to results of federally sponsored research in accordance with the DOE Public Access Plan.


\bibliographystyle{elsarticle-num}
\bibliography{references}



\end{document}